\documentclass[12pt,a4paper,reqno]{amsart}
\usepackage[english]{babel}
\usepackage{amssymb,latexsym,amsfonts,amsthm,upref,amsmath}
\usepackage[foot]{amsaddr}
\usepackage[margin=1in]{geometry} 
\usepackage[dvipsnames,x11names]{xcolor}
 \definecolor{myblue}{HTML}{003399}
\usepackage{hyperref}
\hypersetup{colorlinks,citecolor=myblue,filecolor=black,linkcolor=myblue,urlcolor=myblue}
\usepackage{enumerate}  
\usepackage{tikz}
\usepackage{float}
\usepackage{cleveref}
\usepackage{mathtools}
\usepackage{cite}
\usepackage{filecontents}
\usepackage{enumitem}
\makeatletter
\newcommand{\leqnomode}{\tagsleft@true}
\newcommand{\reqnomode}{\tagsleft@false}

\makeatother
\newtheorem*{thm*}{Theorem}
\newtheorem*{lem*}{Lemma}
\newtheoremstyle{prim}{}{}{\normalfont}{}{\bfseries}{.}{ }{}
\newtheoremstyle{stil}{}{}{\slshape}{}{\bfseries}{.}{ }{}
\theoremstyle{stil}
\newtheorem{thm}{Theorem}[section]
\newtheoremstyle{defi}{}{}{}{}{\bfseries}{.}{ }{}
\theoremstyle{defi}
\newtheorem{defn}[thm]{Definition}
\theoremstyle{defi}
\newtheorem{rem}[thm]{Remark}
\theoremstyle{stil}
\newtheorem*{mthm*}{Main Theorem}
\newtheorem*{kor*}{Corollary}
\newtheorem{pro}[thm]{Proposition}
\theoremstyle{stil}
\newtheorem{lem}[thm]{Lemma}
\theoremstyle{stil}
\newtheorem{kor}[thm]{Corollary}
\theoremstyle{prim}

\newenvironment{prf}{\noindent \textit{Proof.}}{\null\hfill$\qed$\hskip
2mm\vskip 2mm}

\newcommand{\modd}{ \,{\rm mod\,\,}}

\newcommand{\Y}{ {\rm Y}_h^{+}(\g_N)}
\newcommand{\Yht}{ \wht{{\rm Y}}_h^{+}(\g_N)}

\newcommand{\U}{ {\rm U} (\R)}

\newcommand{\V}{\mathcal{V}_{\hspace{-1pt}c}(\g_N)}

\newcommand{\R}{ {\wht{R}}}
\newcommand{\RR}{ {\wtld{R}}}

\newcommand{\vac}{\mathop{\mathrm{\boldsymbol{1}}}}

\newcommand{\gl}{\mathfrak{gl}}

\newcommand{\g}{\mathfrak{g}}

\newcommand{\kp}{\kappa}

\newcommand{\CC}{\mathbb{C}}

\newcommand{\ZZ}{\mathbb{Z}}

\newcommand{\Lc}{\mathcal{L}}

\newcommand{\Sc}{\mathcal{S}}

\newcommand{\Tc}{\mathcal{T}}

\newcommand{\Dc}{\mathcal{D}}

\newcommand{\wtld}{\widetilde}
\newcommand{\wht}{\widehat}

\newcommand{\ot}{\otimes}
\newcommand{\ts}{\hspace{1pt}}

\newcommand{\ndo}{\mathop{\mathrm{End}}}
\newcommand{\om}{\mathop{\mathrm{Hom}}}

\newcommand{\diag}{\mathop{\mathrm{diag}}}
\newcommand{\cdotrl}{\mathop{\hspace{-2pt}\underset{\text{RL}}{\cdot}\hspace{-2pt}}}
\newcommand{\cdotlr}{\mathop{\hspace{-2pt}\underset{\text{LR}}{\cdot}\hspace{-2pt}}}

\newcommand{\jota}{\mathop{\iota_{z_2,z_0,u_1,\ldots ,u_l}}}

\newcommand{\fand}{\quad\text{and}\quad}
\newcommand{\Fand}{\qquad\text{and}\qquad}

\newcommand{\non}{\nonumber}
\newcommand{\beq}{\begin{equation}}
\newcommand{\eeq}{\end{equation}}
\newcommand{\ben}{\begin{equation*}}
\newcommand{\een}{\end{equation*}}

\makeatletter
\def\smalloverbrace#1{\mathop{\vbox{\m@th\ialign{##\crcr\noalign{\kern3\p@}%
  \tiny\downbracefill\crcr\noalign{\kern3\p@\nointerlineskip}%
  $\hfil\displaystyle{#1}\hfil$\crcr}}}\limits}
\makeatother

\makeatletter
\def\smallunderbrace#1{\mathop{\vtop{\m@th\ialign{##\crcr
   $\hfil\displaystyle{#1}\hfil$\crcr
   \noalign{\kern3\p@\nointerlineskip}%
   \tiny\upbracefill\crcr\noalign{\kern3\p@}}}}\limits}
\makeatother

\begin{document}

\title[  $h$-adic quantum vertex algebras   in types $B$, $C$, $D$ and their  $\phi$-coordinated modules]
{$h$-adic quantum vertex algebras   in types $B$, $C$, $D$ and their $\phi$-coordinated modules }

\author{Slaven Ko\v{z}i\'{c}} 
\address{Department of Mathematics, Faculty of Science, University of Zagreb, Bijeni\v{c}ka cesta 30, 10\,000 Zagreb, Croatia}
\email{kslaven@math.hr}

\begin{abstract}
We introduce the $h$-adic quantum vertex algebras  associated with the trigonometric $R$-matrices in types $B$, $C$ and $D$, thus generalizing the well-known Etingof--Kazhdan construction in type $A$. We show that restricted modules for quantum affine algebras in types  $B$, $C$ and $D$ are naturally equipped with the structure of $\phi$-coordinated module for the aforementioned $h$-adic quantum vertex algebras.
\end{abstract}

\maketitle

\allowdisplaybreaks


\section{Introduction}\label{intro}
\numberwithin{equation}{section}
The     vertex algebra  theory presents an important connection between 
 theoretical physics and  mathematics.
It has been extensively studied both by physicists,   in the context  of two-dimensional quantum field theory, and by mathematicians, due to its applications to   representation theory of affine Kac--Moody Lie algebras, finite simple groups and many other areas; see, e.g., the books \cite{FLM,Kac} for   details and references. 
As with vertex algebras, the theory of quantum groups   originates from  
theoretical physics, i.e.,   more specifically, from quantum integrable systems. On the other hand, it possesses  a wide variety of applications to multiple areas of mathematics, such as knot theory, representation theory of algebraic groups in characteristic $p$ etc.; see, e.g., the books \cite{CP,Kas} for   details and references. Motivated by     applications to   two-dimensional statistical models and quantum Yang--Baxter equation, I. Frenkel and Jing \cite{FJ} formulated a fundamental problem of developing the so-called   quantum vertex algebra theory. Roughly speaking, its role has been to associate certain vertex algebra-like objects, i.e. quantum vertex algebras, to various classes of quantum groups, such as quantum affine algebras,   in parallel with the already established connection between affine Kac--Moody
 Lie algebras and vertex algebras.

The notion of quantum vertex algebra, as well as its first examples associated with the rational, trigonometric and elliptic $R$-matrix of type $A$, was introduced by Etingof and Kazhdan  \cite{EK5}. Later on, certain more general structures, such as field algebras and nonlocal vertex algebras, were extensively studied; see, e.g., the papers by Bakalov and Kac \cite{BK} and Li \cite{Li} and references therein. More recently, the structure theory of quantum vertex algebras was further developed by De Sole, Gardini and Kac \cite{DGK} and the corresponding notion of quantum conformal algebra was introduced by  Boyallian and   Meinardi \cite{BM}. 

A major contribution to  the theory   was made by Li \cite{Li1}, who introduced the notion of $\phi$-coordinated modules in order to establish   connection between representation theories of quantum affine algebras and quantum vertex algebras. In comparison with the (usual)
notion of module for (quantum) vertex algebra,  $\phi$-coordinated modules are
characterized by a   deformed version of the weak associativity property which
depends on the choice of the associate $\phi\in\CC((z_2))[[z_0]]$  of the one-dimensional additive formal group. Since then, the theory of $\phi$-coordinated modules has been intensively studied; see, e.g., the recent paper by Jing, Kong, Li and Tan \cite{JKLT} and references therein.

In this paper, we  construct  the $h$-adic quantum vertex algebras  associated with the trigonometric $R$-matrix of type  $B$, $C$, $D$, thus generalizing the aforementioned Etingof--Kazhdan construction in type $A$ \cite{EK5}, along with the corresponding construction for the rational $R$-matrix of type  $B$, $C$, $D$  \cite{BJK}. 
In contrast with the rational $R$-matrix setting, where the quantum vertex algebra structure comes directly from the corresponding double Yangian, the relation between quantum vertex algebras   and quantum affine algebras is no longer so transparent in the trigonometric case.
For example, the original trigonometric $R$-matrix has to be suitably modified, so that, in particular,  it satisfies an additive form  of the Yang--Baxter equation, in order to be used in the aforementioned construction. 
We address such issues in the context   of relations among the   creation and annihilation operators for  the quantum vertex algebras and the   quantum affine algebra generators. 

Regarding the quantum affine algebras, we use their $R$-matrix presentations, as given by Jing, Liu and Molev in types $B$, $C$ and $D$; see \cite{JLM-C,JLM-BD}.
However, we slightly modify the original setting so  that the algebras are defined over the ring $\CC[[h]]$. Following the representation theory of affine Kac--Moody Lie algebras, we consider
a certain wide class of modules for   quantum affine algebras,  the so-called restricted modules. The main result of this paper states that any restricted module for the quantum affine algebra of type $B$, $C$, $D$ is naturally equipped with the structure of $\phi$-coordinated module for the corresponding $h$-adic quantum vertex algebra  with $\phi=z_2 e^{-z_0}$; see Theorem \ref{main1}.
It is proved by a direct  calculation, relying, in  particular, on the $R$-matrix techniques. 
These arguments are easily translated into the rational $R$-matrix setting, thus  proving that any  restricted module  for the double Yangian of type  $A$--$D$ is   equipped with the structure of module  for the $h$-adic quantum vertex algebra associated  with the rational $R$-matrix of the corresponding type, as we explain in Remark \ref{yremark}. 

We should mention that both Theorem \ref{main1} and its converse hold in the case of the trigonometric $R$-matrix in type $A$ with $\phi=z_2 e^{z_0}$; see \cite{Koz1}. However, their proofs  rely on   Ding's quantum current realization \cite{D}, which is known only in type $A$. At the end of the paper, we discuss quantum current commutation relations in types $B$, $C$ and $D$, which might   lead to the  appropriate definition of quantum current algebras in these types and, hopefully, to the proof of the converse of Theorem \ref{main1}.

\section{Trigonometric $R$-matrix}\label{sec02}
%

In this section, we consider   the  trigonometric $R$-matrix of types $B$, $C$ and $D$. Our main focus is on its   additive  form, defined over the ring $\CC[[h]]$, which    satisfies the additive versions of  the Yang--Baxter equation, crossing symmetry and unitarity properties; see Proposition \ref{Rpro}. 
We start by introducing some notation, which mostly follows \cite[Sect. 1]{JLM-C} and \cite[Sect. 1]{JLM-BD}. Let $\g_N$  be  the orthogonal Lie algebra  $\mathfrak{o}_N$   or the symplectic Lie algebra $\mathfrak{sp}_N$, where $N$ is even in the symplectic case. Write $n=\lfloor{N/2\rfloor}$ so that $\g_N$ is of type $B_n, C_n,D_n$, i.e. so that   we have
$\g_N =\mathfrak{o}_{2n+1}, \mathfrak{sp}_{2n}, \mathfrak{o}_{2n}$,  respectively. 
For $i=1,\ldots ,N$ let $i'=N+1-i$. 
In the symplectic case, we introduce the sign
$$\varepsilon_i=\begin{cases}
1,&\text{if }i=1,\ldots ,n,\\ 
-1,&\text{if }i=n+1,\ldots ,2n.
\end{cases}
$$
To consider all cases simultaneously 
we  set $\varepsilon_i=1$ for   $i=1,\ldots ,N$ in the orthogonal case. 
Define the matrix transposition on $\ndo\CC^N$ by 
$e_{ij}^t = \varepsilon_i\varepsilon_j e_{j' i'}$, where $e_{ij} $ denote matrix units. Finally, introduce the $N$-tuple
\beq\label{frac6}
(\bar{1},\ldots , \bar{N})=\begin{cases}
(\textstyle n-\frac{1}{2},\ldots,\frac{3}{2},\frac{1}{2},0,-\frac{1}{2},-\frac{3}{2},\ldots ,-n+\frac{1}{2}),&\text{if }\g_N =\mathfrak{o}_{2n+1},\\ 
(n,\ldots,2,1,-1,-2,\ldots ,-n),&\text{if }\g_N =\mathfrak{sp}_{2n},\\
(n-1,\ldots , 2,1,0,0,-1,-2,\ldots,-n+1),&\text{if }\g_N =\mathfrak{o}_{2n}.
\end{cases}
\eeq

Let $q$ be a formal parameter and
\beq\label{xi}
\xi =\begin{cases}
q^{2-N},&\text{if }\g_N =\mathfrak{o}_{N},\\
q^{-2-N},&\text{if }\g_N =\mathfrak{sp}_{N}.
\end{cases}
\eeq
As in \cite{FR}, let   
\beq\label{freg}
f(x)= 1+\sum_{k= 1}^\infty f_k x^k,\quad\text{where}\quad f_k\in\CC(q)\text{ are regular at }q=1 , 
\eeq
be
a unique formal power series in $ \CC(q)[[x]] $
such that
\beq\label{uk3}
f(x)f(x\xi)=\frac{1}{\left(1-xq^{-2}\right)\left(1-xq^2\right)\left(1-x\xi\right)\left(1-x\xi^{-1}\right)}.
\eeq
In accordance with \cite{J}, define the $R$-matrix 
\begin{gather}
R(x,q)=f(x)\ts 
R^+ (x,q),
\qquad
\text{where}\label{R}\\
R^+ (x,q) =
q^{-1}(x-1)(x-\xi)R 
-(q^{-2}-1)(x-\xi)P
+\xi (q^{-2}-1)(x-1)  Q\label{erplus}
\end{gather}
and the operators $P$, $Q$ and $R$   are given by
\begin{gather*}
P=\sum_{i,j=1}^N e_{ij}\ot e_{ji},
\qquad
Q=\sum_{i,j=1}^N q^{\bar{i}-\bar{j}}\varepsilon_i\varepsilon_j e_{i'j'}\ot e_{ij}\Fand\\
R= \,\,
q\sum_{i=1,\ldots ,N , \,i\neq i'} e_{ii}\ot e_{ii}
+\sum_{i\neq j,j'} e_{ii}\ot e_{jj}
+q^{-1}\sum_{i\neq i'} e_{ii}\ot e_{i'i'}+(q-q^{-1})\sum_{i<j}e_{ij}\ot e_{ji}\\
 \hspace{-65pt}-(q-q^{-1})\sum_{i>j}q^{\bar{i}-\bar{j}}\varepsilon_i\varepsilon_j e_{i'j'}\ot e_{ij}
+\delta_{N\ts 2n+1}e_{n+1\ts n+1}\ot e_{n+1\ts n+1}.
\end{gather*}

Throughout this paper we use the standard tensor notation, where for any
$$
A=\sum_{i,j,k,l=1}^N a_{ijkl}\ts e_{ij}\ot e_{kl}
$$
and indices $1\leqslant r,s\leqslant m$ such that $r\neq s$, we write $A_{rs}$ for the element
\beq\label{nottat}
A_{rs}=\sum_{i,j,k,l=1}^N a_{ijkl}\ts (1^{\ot(r-1)}\ot e_{ij}\ot 1^{\ot(m-r)}) (1^{\ot(s-1)}\ot e_{kl}\ot 1^{\ot(m-s)}).
\eeq
To simplify the notation we suppress the dependence on the parameter $q$ 
and write
$R(x)$ instead of $R(x,q) $.
The $R$-matrix \eqref{R} satisfies the Yang--Baxter equation
\beq\label{YBE}
R_{12}(x) R_{13}(xy) R_{23}(y)
=R_{23}(y) R_{13}(xy)R_{12}(x).
\eeq
Here we used the notation convention   \eqref{nottat} with $m=3$ and $(r,s)\in\left\{(1,2),(1,3),(2,3)\right\}$. Furthermore, $R(x)$ possesses the
  crossing symmetry property 
\beq\label{CS}
R(x) D_1 R(x\xi)^{t_1} D_1^{-1} =\xi^2 q^{-2},
\eeq
where    $D=\diag(q^{\bar{1}},\ldots, q^{\bar{N}})$ is the diagonal matrix and the subscript $1$ indicates that the matrix $D$ and the transposition ${}^t$ are applied on the first tensor factor of
$\ndo\CC^N \ot\ndo\CC^N$.

Clearly, the entries of the $R$-matrix $R(x)$  are formal power series in $\CC(q^{1/2})[[x]]$, where the fractional powers of $q$ appear only when $\g_N=\mathfrak{o}_{2n+1}$; recall \eqref{frac6}. We  now want to apply the substitutions
\beq\label{sub}
x=e^{-u}=\sum_{k\geqslant 0}  \frac{(-u)^k}{k!}\in\CC[[u]] \fand q^{1/2}=e^{h/4}=\sum_{k\geqslant 0}\frac{h^k}{4^k k!}\in\CC[[h]]
\eeq
 to  \eqref{R}, where $u$ is a variable and $h$ a formal parameter, so that we obtain the $R$-matrix with entries in $\CC((u))[[h]]$. 
Let $\kp=\frac{N}{2}\pm 1$, where the plus (minus) sign appears in the symplectic (orthogonal) case,
so that we have  $e^{-\kp h} =\xi|_{q=e^{h/2}}$. 

\begin{rem}\label{rem_exp}
Note that  by applying substitutions \eqref{sub} to the expressions of the form $(q^a -x  q^b)^{r}$ with $r<0$ and then using  the usual   expansions
we obtain elements of $\CC((u))[[h]]$. Indeed, for any $a,b\in\frac{1}{2}\ZZ$ and a negative integer $r$ we have
\begin{align}
&(q^a -x  q^b)^{r}\big|_{x=e^{-u},q^{1/2}=e^{h/4}}
= e^{ ar h/2}
(-u+(b-a)h/2)^{r} G(u,h), \qquad\text{where}\label{tmp1}\\
&\qquad (-u+(b-a)h/2)^{r} 
= \sum_{k\geqslant 0} \binom{r}{k} (-u)^{r-k} \left(\frac{(b-a)h}{2}\right)^k
\in\CC[u^{-1}][[h]],\non\\
&\qquad G(u,h)=\left(\frac{ -u+(b-a)h/2  }{1-e^{-u+(b-a)h/2}}\right)^{-r}\in\CC[[u,h]],\non
\end{align}
so that \eqref{tmp1} belongs to $\CC((u))[[h]]$.  
\end{rem}

Let $F_{*}(x_1,\ldots ,x_l)$ be the localization of the ring of Taylor series $F[[x_1,\ldots ,x_l]]$ over the field $F$ at  nonzero polynomials 
$F[x_1,\ldots ,x_l]^{\times}$. Denote by $\iota_{x_1,\ldots ,x_l}$   the unique embedding  
$F_{*}(x_1,\ldots ,x_l)\to F((x_1))\ldots((x_l))$. 
 Note that in \eqref{tmp1} we  implicitly  applied the embedding $\iota_{u,h}$ with $F=\CC$, thus getting an element of $\CC((u))[[h]]\subset\CC((u))((h))$. Also, the right hand side of \eqref{uk3} was understood as an element of $\CC(q)[[x]]$ via the embedding $\iota_x$, where $F=\CC(q)$.
As with \eqref{uk3} and \eqref{tmp1}, throughout the  paper  we shall often omit the embedding symbol and employ the  usual  expansion convention where the   embedding is determined by the order of the variables: if $\pi $ is a permutation in the symmetric group $\mathfrak{S}_l$, then
$(x_{\pi_1}+\ldots +x_{\pi_l})^{r}$ with $r<0$ stands for $\iota_{x_{\pi_1},\ldots ,x_{\pi_l}}(x_{\pi_1}+\ldots +x_{\pi_l})^{r}$.

\begin{lem}\label{l1}
There exists a unique formal power series
$g(u,h)\in  \CC((u))[[h]]$
such that $g (u,0)= (1-e^{-u})^{-2}$ and such that
\beq\label{tmp2}
g(u,h)g(u+\kp h,h)=\frac{1}{(1-e^{-u-h})(1-e^{-u+h})(1-e^{-u-\kp h})(1-e^{-u+\kp h})}.
\eeq
\end{lem}

\begin{prf}
Let   $
g(u,h)=\sum_{l\geqslant 0} g_l(u)h^l$ for some $g_{l}(u)\in\CC((u))
$. By using the formal Taylor theorem we   write the left hand side of  \eqref{tmp2} as
$$
\left(\sum_{l\geqslant 0} g_l(u)h^l\right)
\left(\sum_{k\geqslant 0} \left(\sum_{m=0}^k \frac{ \kp^m  }{m!}\frac{d^m}{du^m}g_{k-m} (u)\right)h^k  \right).
$$
On the other hand, by employing \eqref{tmp1}, the right hand side of  \eqref{tmp2} is expressed as a  power series in   $\CC((u))[[h]]$. The series $g_l(u)$ are now calculated by comparing
  the coefficients of $h^l$ for $l=0,1, \ldots$ on the left and the right hand side of the equality.
\end{prf}

\begin{rem}\label{rem12}
By   computing the coefficients $g_l(u)$
from the identity established in the proof of Lemma \ref{l1}, one finds that  the series $g(u,h)$ admits the form
\beq\label{guh1}
g(u,h)=\sum_{l\geqslant 0} \frac{p_l (e^{-u})}{ (1-e^{-u})^{r_l}}h^l \quad\text{for some}\quad p_l(z) \in\CC[z],\,r_l\in\ZZ_{\geqslant 0}.
\eeq
Indeed, this is   verified by induction over $l$ which    relies on the Taylor  expansions
$$
\frac{1}{1-e^{-(u+a h)}}
=\sum_{l\geqslant 0}\frac{(a h)^l}{l!} \frac{\partial^l}{\partial u^l} \frac{1}{1-e^{-u}}, \quad\text{where }a\in\textstyle\frac{1}{2}\ZZ.
$$
Hence, the substitution $e^{-u}=z$ in \eqref{guh1} produces a series
$g_1(z,h)\in  \CC[[z,h]]$ satisfying
\beq\label{uk4}
g_1 (z,h)g_1(ze^{-\kp h},h) 
=\frac{1}{(1-ze^{-h})(1-ze^{h})(1-ze^{-\kp h})(1- ze^{\kp h})}.
\eeq
Furthermore, the condition $g (u,0)= (1-e^{-u})^{-2}$, along with \eqref{uk4} at $z=0$, implies that $ g_1(z,h)$ belongs to $ 1+z\CC[[z,h]]$. However, as with \eqref{uk3}, one easily checks that the equation \eqref{uk4} possesses a unique solution in $ 1+z\CC[[z,h]]$, so we conclude that
\beq\label{rem12eq}
g_1 (z,h) = f(z)\big|_{q=e^{h/2}}.\big.
\eeq
Therefore, by \eqref{guh1} the    power series $f(z)\big|_{q=e^{h/2}} \big.$ admits the form
\beq\label{uk6}
f(z)\big|_{q=e^{h/2}} \big. =\sum_{l\geqslant 0} \frac{p_l (z)}{(1-z)^{r_l}}h^l \quad\text{for some}\quad p_l(z) \in\CC[z],\,r_l\in\ZZ_{\geqslant 0}.
\eeq
Finally,  by applying the substitution $z=e^{-u}$ to \eqref{uk6}  we find
\beq\label{twolines}
f(z)\big|_{q=e^{h/2}, z=e^{-u}} \big. =g(u,h).
\eeq
\end{rem}

We single out the following simple consequence of \eqref{uk6} which will be useful later on.
\begin{kor}\label{korf}
For any integer $k> 0$ there exists an integer $r\geqslant 0$ such that
$$
(1-z)^r  f(z)\big|_{q=e^{h/2}}\in\CC[z,h]\mod h^k.
$$
\end{kor}
In other words, the corollary asserts that   the coefficients of the powers $h^0,h^1,\ldots,h^{k-1}$ of the given expression belong to $\CC[z]$. Clearly, this holds for any $r\geqslant\max\left\{r_0,\ldots,r_{k-1}\right\}$.

Consider the $R$-matrix 
\begin{align}
\R(u)=&\,\,
e^{(1+2\kp)h/2}\ts g(u,h) \ts R^+(e^{-u},e^{h/2}).\label{Rh}
\end{align}
The above formula, up to the multiplicative factor $e^{(1+2\kp)h/2} $, is   obtained  by applying the substitutions \eqref{sub} to the $R$-matrix $R(x,q)$   defined by  \eqref{R} in the sense of \eqref{uk6} and  \eqref{twolines}.
 Lemma  \ref{l1} implies that $\R(u)$ belongs to $\ndo\CC^N\ot \ndo\CC^N ((u))[[h]]$. Its main properties are listed in the following proposition.

\begin{pro}\label{Rpro}
The $R$-matrix \eqref{Rh} satisfies the Yang--Baxter equation 
\beq\label{YBEh}
\R_{12}( u ) \R_{13}(u+v) \R_{23}(v)
=\R_{23}(v) \R_{13}(u+v)\R_{12}(u),
\eeq
the crossing symmetry relation     
\beq\label{CSh}
\R( u ) M_1 \R( u+\kp h )^{t_1} M_1^{-1} = 1 ,\quad\text{where}\quad M=D|_{q^{1/2}=e^{h/4}} ,
\eeq
and the unitarity relation
\beq\label{Uh}
\R_{12}(u) \R_{21}(-u)=1.
\eeq
\end{pro}

\begin{prf}
The Yang--Baxter equation \eqref{YBEh} is a direct consequence of \eqref{YBE}. The proof of \eqref{CSh} and \eqref{Uh}   relies on the properties of the $R$-matrix $\bar{R}(x)$, defined by the identity
$$
R (x )=f(x)(x-q^{-2})(x-\xi)\bar{R}(x).
$$
More specifically,  by \cite[Sect. 3.1]{JLM-C} and \cite[Sect. 3.1]{JLM-BD}    the aforementioned $R$-matrix possesses the following crossing symmetry and unitarity properties,
\beq\label{CSbar}
\bar{R}(x) D_1 \bar{R}(x\xi)^{t_1} D_1^{-1} = \frac{(x-q^2)(x\xi -1)}{(1-x)(1-x\xi q^2)}\Fand \bar{R}_{12}(x) \bar{R}_{21}(x^{-1})=1.
\eeq

Let us prove  \eqref{CSh}. By the crossing  symmetry property in \eqref{CSbar} it is sufficient to check that the series $g(u,h)$, as given by Lemma \ref{l1}, satisfies the identity
\begin{align*}
&  g(u,h)g(u+\kp h,h)(e^{-u}-e^{-h})(e^{-u} -e^{-\kp h})(e^{-u-\kp h}-e^{-h})(e^{-u-\kp h} -e^{-\kp h})\\
 &\qquad= e^{(-2\kp-1)h}\frac{(1-e^{-u})(1-e^{-u+(-\kp +1)h})}{(e^{-u}-e^h)(e^{-u-\kp h}-1)}.
\end{align*}
However, this follows directly from \eqref{tmp2}.

As for  \eqref{Uh}, by the unitarity property in \eqref{CSbar} it is sufficient to show that
\beq\label{tmp4}
e^{(1+2\kp)h } g(u,h)g(-u,h)(e^{-u}-e^{-h})(e^{-u} -e^{-\kp h})(e^u-e^{-h})(e^u -e^{-\kp h})=1.
\eeq
Denote the left hand side of \eqref{tmp4} by $G(u,h)$. Clearly, the power series $G(u,h)$ belongs to $\CC((u))[[h]]$. By using \eqref{tmp2} one easily checks the identity
\beq\label{tmp5}
G(u,h)G(u+\kp h, h)=1.
\eeq
As with the proof of Lemma \ref{l1}, write
$ 
G(u,h)=\sum_{k\geqslant 0} G_k(u) h^k$ with $G_k(u)\in\CC((u))$, 
so that the application of    the formal Taylor
theorem  turns \eqref{tmp5} into
$$
\left(\sum_{k\geqslant 0} G_k(u)h^k\right)
\left(\sum_{l\geqslant 0} \left(\sum_{m=0}^l \frac{ \kp^m  }{m!}\frac{d^m}{dx^m}G_{l-m} (u)\right)h^l  \right)=1.
$$
By comparing the coefficients of $h^n$, $n\geqslant 0$, on the left and the right hand side we find that $G_k(u)=\pm \delta_{k0}$, so that $G(u,h)=\pm 1$. Finally, by Lemma \ref{l1} we have
$g (u,0)= (1-e^{-u})^{-2}$ so it is clear from   \eqref{tmp4}  that $G(u,0)=  1$. Thus, we have    $G(u,h)=  1$, as required.
\end{prf}

We finish this section by recalling some  notation from \cite[Sect. 2.2]{Koz1}.
Suppose $V$ is a topologically free $\CC[[h]]$-module and $a_1,\ldots ,a_l,k >0$ are    integers. Let $A$ be an   element of $X\coloneqq \om(V,V[[z_1^{\pm 1},z_2^{\pm 1},u_1,\ldots ,u_l]])$ such that there exist elements
$B$  in $\om(V,V((z_1,z_2))[[u_1,\ldots ,u_l,h]])$ and
 $C_1,\ldots , C_{l+1}\in X$  so that we have
\begin{gather}
A=B+ \sum_{i=1}^l u_i^{a_i} C_i +  h^k C_{l+1}.\label{decompmod}
\end{gather}
To indicate   that $A $ admits such decomposition, we  shall write
$$
A\in\om(V,V((z_1,z_2))[[u_1,\ldots ,u_l]] )\mod u_1^{a_1},\ldots, u_l^{a_l},h^k.
$$
Later on, we shall use the fact that, e.g.,  the substitution
\beq\label{bsupst}
B\big|_{z_1= z_2 e^{-z_0} }\big. =\jota\left( B(z_1,z_2,u_1,\ldots , u_l)\big|_{z_1= z_2 e^{-z_0}}\big.\right)
\eeq
is well-defined even though the  substitution  $\textstyle A\big|_{z_1= z_2 e^{-z_0}}\big.$ does not exist in general. 
Note that the element $B\in\om(V,V((z_1,z_2))[[u_1,\ldots ,u_l,h]])$ satisfying \eqref{decompmod}  is    unique modulo
$ 
\sum_{i=1}^l u_i^{a_i}X + h^k X 
$. 
To   simplify the notation, we shall denote \eqref{bsupst} as
$$
A \big|_{z_1= z_2 e^{-z_0}}^{\modd u_1^{a_1},\ldots, u_l^{a_l},h^k}\big. =A(z_1,z_2,u_1,\ldots , u_l)\big|_{z_1= z_2 e^{-z_0}}^{\modd u_1^{a_1},\ldots, u_l^{a_l},h^k}\big. .
$$

\section{Quantum affine algebras in types $B$, $C$ and $D$}\label{sec05}
%

In this section, we  follow the exposition in \cite[Sect. 1]{JLM-C} and \cite[Sect. 1]{JLM-BD} to recall the $R$-matrix presentations  of   quantum affine algebras in types $B$, $C$ and $D$.
 However, we make some minor adjustments  so that, in contrast with the original presentations, the quantum affine algebra structure is defined over the commutative ring $\CC[[h]]$  instead over the field $\CC(q^{1/2})$.
At the end, we introduce   the notion of restricted module  for   quantum affine algebra   in parallel with the  representation theory of  affine Kac--Moody Lie algebras.

The quantum affine algebra defining relations are expressed using the   $R$-matrix 
\beq\label{rrr1}
\RR(x)=e^{( 2\kp +1)h/2} \ts R(x,q)\left|_{q^{1/2}=e^{h/4}},\right.
\eeq
where $R(x,q)$ is defined by \eqref{R}. Note that  $\RR(x)$ is a well-defined element of 
$\ndo\CC^N \ot\ndo\CC^N [[x,h]]$ due to discussion in Remark \ref{rem12}; recall, in particular, \eqref{rem12eq}.
 Moreover, the   $R$-matrix \eqref{rrr1}  satisfies the Yang--Baxter equation, 
\beq\label{YBEv2}
\RR_{12}(x) \RR_{13}(xy) \RR_{23}(y)
=\RR_{23}(y) \RR_{13}(xy)\RR_{12}(x) 
\eeq
and the crossing symmetry relation,
\beq\label{CSv2}
\RR(x) M_1 \RR(x\zeta)^{t_1} M_1^{-1} =1,\quad\text{where}\quad \zeta = \xi\left|_{q=e^{h/2}}\right. ,
\eeq
 which follows from   \eqref{YBE} and \eqref{CS}, respectively. Note that, comparing to \eqref{CS}, the crossing symmetry property in \eqref{CSv2} takes the slightly simpler form due to the 
normalization term $e^{( 2\kp +1)h/2}$ in \eqref{rrr1}. Clearly, such normalization does not affect the quantum affine algebra defining relations below. However, it establishes the following direct correspondence between the $R$-matrices \eqref{Rh} and \eqref{rrr1}:

\begin{pro}\label{imp_prop}
For any integers $a,b,l>0$ and $\alpha\in\CC$ there exists an integer $r\geqslant 0$ such that the coefficients of all monomials
\beq\label{monomi}
u^{a'} v^{b'} h^{l'},\qquad\text{where}\qquad 0\leqslant a'<a,\quad 0\leqslant b'<b \fand 0\leqslant l'<l,
\eeq
in $(x-y)^r \RR(xe^{u-v+\alpha h} /y)$ belong to $\ndo\CC^N\ot\ndo\CC^N [x,y^{\pm 1}]$ and such that the coefficients of all monomials   \eqref{monomi} in
$$
\left((x-y)^r \RR(xe^{u-v+\alpha h} /y)\right)\Big|_{y=xe^{-z}}^{\text{mod }u^a, v^b, h^l}
\fand
x^r \left(1-e^{-z}\right)^r\R(-z-u+v-\alpha h)
\Big.$$
coincide.
\end{pro}

\begin{prf}
Clearly,   the statement of the proposition holds if the $R$-matrices 
$ \RR(xe^{u-v+\alpha h} /y)$ and $\R(-z-u+v-\alpha h)$ 
are replaced by
$ 
R^+(xe^{u-v+\alpha h} /y,e^{h/2})$ and 
$R^+(e^{z+u-v+\alpha h},e^{h/2})$,
respectively; recall \eqref{erplus}. Therefore, it is sufficient to verify the corresponding statement for the normalizing series $g_1 (z,h) = f(z)\big|_{q=e^{h/2}}$ and $g(u,h)$ of the   $R$-matrices   \eqref{rrr1} and \eqref{Rh}, which are given by \eqref{rem12eq} and Lemma \ref{l1}, respectively. However, this follows by an argument, relying on Remark \ref{rem12} and Corollary \ref{korf}, which goes in parallel with the proof of  \cite[Lemma 3.2]{Koz1}, where the same property for the normalizing series of the trigonometric $R$-matrix in type $A$ was established. 
\end{prf}

The {\em quantum affine algebra $U_h(\wht{\g}_N)_c$ at the level $c$}, where $c\in\CC$ and $\g_N $ is the Lie algebra $\mathfrak{o}_{2n+1}, \mathfrak{sp}_{2n}, \mathfrak{o}_{2n}$, is defined as the $h$-adically complete  associative algebra over the ring $\CC[[h]]$ generated by  the elements $l_{ij}^{\pm }(\mp r)$ with $i,j=1,\ldots ,N$ and $r=0,1,\ldots $ such that for all $i=1,\ldots ,N$
$$
l_{ij}^+(0)=l_{ji}^-(0)=0\quad\text{for}\quad i<j\Fand
 l_{ii}^-(0)-l_{ii}^+(0)=hl_{ii}^+(0)l_{ii}^-(0)=hl_{ii}^-(0)l_{ii}^+(0).
$$
 The generators are subject to  defining relations which are expressed  as follows. Let
$$
L^{\pm }(x)=\sum_{i,j=1}^N e_{ij}\ot l_{ij}^{\pm}(x),\quad\text{where}\quad
l_{ij}^{\pm}(x)=\delta_{ij}\mp h\sum_{r=0}^\infty l_{ij}^{\pm}(\mp r) x^{\pm r}. 
$$
The defining relations are given by
\begin{gather}
\RR(x/y)\ts L_1^{\pm}(x)\ts L_2^{\pm }(y) =L_2^{\pm }(y)\ts L_1^{\pm}(x)\ts \RR(x/y),\label{uqrtt1}\\
\RR(xe^{hc/2}/y)\ts L_1^{+}(x)\ts L_2^{- }(y) =L_2^{- }(y)\ts L_1^{+}(x)\ts \RR(xe^{-hc/2}/y),\label{uqrtt2}\\
L^{\pm}(x)\ts M \ts L^{\pm}(x\zeta)^t\ts M^{-1} =1,\label{uqrtt3}
\end{gather}
where the transposition is defined as in Section \ref{sec02}, the diagonal matrix $M $ is given by \eqref{CSh} and $\zeta = \xi\left|_{q=e^{h/2}}\right.$; recall \eqref{xi}. The subscripts $1$ and $2$ in \eqref{uqrtt1} and \eqref{uqrtt2} indicate the corresponding tensor copies of $\ndo\CC^N$ in $\ndo\CC^N\ot \ndo\CC^N\ot U_h(\wht{\g}_N)_c$ so that
$$
L_1^{\pm }(x)=\sum_{i,j=1}^N e_{ij}\ot 1\ot l_{ij}^{\pm}(x)
\fand
L_2^{\pm }(x)=\sum_{i,j=1}^N 1\ot e_{ij}\ot  l_{ij}^{\pm}(x).
$$
More generally, in accordance with \eqref{nottat}, for any indices $1\leqslant r\leqslant m$ we shall write
\beq\label{not227}
L_r^{\pm }(x)=\sum_{i,j=1}^N 1^{\ot (r-1)}\ot e_{ij}\ot 1^{\ot(m-r)}\ot l_{ij}^{\pm}(x).
\eeq

It is worth noting that the defining relations \eqref{uqrtt1} and  \eqref{uqrtt3} are equivalent to 
\begin{gather}
  R^+(x/y,e^{h/2})\ts L_1^{\pm}(x)\ts L_2^{\pm }(y) =L_2^{\pm }(y)\ts L_1^{\pm}(x)\ts R^+(x/y,e^{h/2}),\label{equqrtt1}\\
 L^{\pm}(x\zeta )  \cdotlr     \left(M\ts L^{\pm}(x)^t  \right)   =M
\label{equqrtt3},
\end{gather}
respectively. The subscript LR in \eqref{equqrtt3} indicates that the first (second) tensor factor of $ L^{\pm}(x\zeta ) $ is applied from the left (right) to $\left(M\ts L^{\pm}(x)^t  \right)$. We shall often use such ordered product notation. In addition, \eqref{equqrtt3} can be equivalently written as
\beq\label{nmh}
\left(M\ts L^{\pm}(x)^t  \right)
\cdotrl   
L^{\pm}(x\zeta )    =M,
\eeq
where the
subscript RL now indicates that the first (second) tensor factor of $\left(M\ts L^{\pm}(x)^t  \right)$ is applied from the right (left). Clearly, both \eqref{equqrtt3} and \eqref{nmh} are found  by  multiplying \eqref{uqrtt3} by $M$ from the right and then applying the transposition.
 
The $RLL$-relations \eqref{uqrtt1} and \eqref{uqrtt2} can be       generalized
 as follows.
First, for any  $a\in\CC$, integers $m,k\geqslant 1$ and the  variables $x=(x_1,\ldots ,x_k)$, $y=(y_1,\ldots ,y_m)$ define
\beq\label{formulla3}
\RR_{km}^{12}( xe^{ah}/y )= \prod_{i=1,\dots,k}^{\longrightarrow} 
\prod_{j=k+1,\ldots,k+m}^{\longleftarrow} \RR_{ij}( x_i e^{ah}/y_{j-k}) 
\eeq
with the arrows indicating the order of the factors. 
Using   \eqref{not227} introduce the elements
$$
L_{[k]}^{\pm}(x )  = L_1^\pm (x_1)\ldots L_k^\pm (x_k)\in(\ndo\CC^N)^{\ot k}\ot U_h(\wht{\g}_N)_c [[x_1^{\pm 1},\ldots ,x_k^{\pm 1}]].
$$
Defining relations   \eqref{uqrtt1} and \eqref{uqrtt2} imply the identities
\begin{gather}
\RR_{km}^{12}(x/y)\ts L_{[k]}^{\pm 13}(x)\ts L_{[m]}^{\pm 23}(y) =L_{[m]}^{\pm 23}(y)\ts L_{[k]}^{\pm 13}(x)\ts \RR_{km}^{12}(x/y),\label{guqrtt1}\\
\RR_{km}^{12}(xe^{hc/2}/y) \ts L_{[k]}^{+13}(x)\ts L_{[m]}^{-23 }(y) =L_{[m]}^{-23 }(y)\ts L_{[k]}^{+13}(x)\ts \RR_{km}^{12}(xe^{-hc/2}/y),\label{guqrtt2} 
\end{gather}
where the superscripts $1,2,3$   indicate the tensor factors as follows:
$$
\smalloverbrace{(\ndo\CC^N)^{\ot k}}^{1} \ot \smalloverbrace{(\ndo\CC^N)^{\ot m}}^{2}\ot \smalloverbrace{U_h(\wht{\g}_N)_c}^{3}.
$$

 Let $W$ be a $\CC[[h]]$-module. We denote  by $W((x))_h$ the $\CC[[h]]$-module of all series  
\beq\label{notacija}
a(x)=\sum_{r\in\ZZ } a_r x^{-r-1}\in W[[x^{\pm 1}]]\quad \text{such that}\quad a_r \rightarrow 0 \text{ when } r \rightarrow \infty
\eeq
 with respect to the $h$-adic topology. 
Moreover, let $W[x^{-1}]_h$ be the $\CC[[h]]$-module of all power series as in \eqref{notacija} such that, in addition, $a(x)$ belongs to $W[[x^{-1}]]$.
This notation naturally extends to the multiple variable case,   so   we  write, for example, $W((x_1,\ldots ,x_n))_h $.

In this paper, we shall consider the so-called restricted modules for quantum affine algebras. An $U_h(\wht{\g}_N)_c$-module $W$ is said to be {\em restricted} if it is topologically free as a $\CC[[h]]$-module and, in addition,  the action of $L^-(x)$ on $W$ satisfies
\beq\label{restricted1}
L^-(x) \in \ndo\CC^N \ot \om(W,W[x^{-1}]_h).
\eeq
Requirement  \eqref{restricted1} can be equivalently expressed as follows: for any $w\in W$ and integer $k>0$ the expression $L^-(x)w$ contains only finitely many powers of $x^{-1}$ modulo $h^k$. 
More generally, if $W$ is restricted, we have
\beq\label{restr5}
L_{[k]}^{-}(x_1,\ldots ,x_k)  \in(\ndo\CC^N)^{\ot k}\ot \om(W,W [x_1^{- 1},\ldots ,x_k^{- 1}]_h).
\eeq 
Finally, it is worth noting that, even though the products such as $L^+(x)L^-(x)$ are not  defined in the quantum affine algebra, they can be regarded as   well-defined elements of $\ndo\CC^N \ot\om(W,W((x))_h)$ if $W$ is  a restricted module, which we shall use later on.

\section{$h$-adic quantum vertex algebras in types $B$, $C$ and $D$}\label{sec03}

\subsection{Constructing  $h$-adic quantum vertex algebras}\label{subsec0304}
In this subsection we introduce the $h$-adic quantum vertex algebras associated with the trigonometric $R$-matrices of types $B$, $C$ and $D$.

\subsubsection{Creation operators}\label{subsec0301}

We   follow the approach from \cite[Sect. 3]{EK3} to associate an algebra to the $R$-matrix \eqref{Rh}.
The quantized universal enveloping algebra $\U $ is defined as the
topologically free  associative algebra over the ring $\CC[[h]]$ generated by the elements $t_{ij}^{(-r)}$, where $i,j=1,\ldots ,N$ and $r=1,2,\ldots ,$ subject to the defining relations
\begin{gather}
\R(u-v)\ts T_1^+(u)\ts T_2^+(v)=T_2^+(v)\ts T_1^+(u)\ts \R(u-v),\label{RTTh}\\
T^+(u)\ts  M\ts  T^+(u+\kp h)^t\ts  M^{-1} =1.\label{RELh}
\end{gather} 
Here $T^+(u)$ denotes the matrix  
\beq\label{tplusu}
T^+(u)=\sum_{i,j=1}^N e_{ij}\ot t_{ij}^+(u) ,\quad\text{where}\quad
t_{ij}^+(u)=\delta_{ij}- h\sum_{s\geqslant 1} t_{ij}^{(-s)} u^{s-1}.
\eeq
Clearly, the defining relations \eqref{RTTh} are equivalent to 
\beq\label{eqRTTh}
 R^+(e^{-u+v},e^{h/2})\ts T_1^+(u)\ts T_2^+(v)=T_2^+(v)\ts T_1^+(u) \ts  R^+(e^{-u+v},e^{h/2}).
\eeq

\begin{rem}\label{tfree-remark}
Recall that  $\U $ is defined as a topologically free algebra over $\CC[[h]]$; cf. \cite[Ch. XVI]{Kas}. Suppose $F$ is the free $\CC[[h]]$-algebra in the given generators and $I$ is the   ideal of the defining relations in $F$.  Write  $[I]$ for the ideal of all $x\in F$ such that $h^m x$ belongs to $I$ for some integer $m\geqslant 0$. The algebra $\U $ is then defined as the quotient of the $h$-adic completion of $F$ by the $h$-adic completion of $[I]$. By arguing as in the proof of \cite[Prop. 2.2]{c11} one checks that $\U $ is topologically free. 
\end{rem}

\begin{rem}\label{remark31}
In this remark we discuss   connection between the defining relations for the quantized universal enveloping algebra and the quantum affine algebra structure.
Denote by $U'_h(\wht{\g}_N)_c$ the quotient of the quantum affine algebra $U_h(\wht{\g}_N)_c$ over its $h$-adically closed ideal generated by the elements $l_{ii}^+(0)$, $i=1,\ldots ,N$.
Let
$\Y$
be the $h$-adically completed subalgebra of $U'_h(\wht{\g}_N)_c$ generated by the elements $l_{ij}^+(-r)$ with $i,j=1, \ldots ,N$ and $r=0,1,\ldots . $ Here we use the same notation for the elements of     $U_h(\wht{\g}_N)_c$ and for their images in the quotient $U'_h(\wht{\g}_N)_c$. As in Section \ref{sec05}, we organize the generators of the algebra $\Y$ into matrices of formal power series
$$
L^{+ }(x)=\sum_{i,j=1}^N e_{ij}\ot l_{ij}^{+}(x),\quad\text{where}\quad
l_{ij}^{+}(x)=\delta_{ij}- h\sum_{r=0}^\infty l_{ij}^{+}(- r) \ts x^{ r} . 
$$
Note that   the constant term of $L^+(x)$  is a lower triangular matrix with units on its main diagonal since $l_{ij}^+(0)=0$ for $i\leqslant j$ in $\Y$.
Denote by $\Yht$ the completion of the algebra $\Y$ with respect to  the descending filtration defined by setting the degree of the generators $l_{ij}^{+}(- r)$ to be equal to $r$.
Consider the matrix 
$$
\Tc^+(u)=\sum_{i,j=1}^N e_{ij}\ot\tau^+_{ij}(u),\quad\text{where}\quad
\tau^+_{ij}(u)=\delta_{ij}-h\sum_{s \geqslant 1} \tau^{(-s)}_{ ij }u^{s-1},
$$
defined by
$
\Tc^+(u)=L^+(e^{-u}) .
$
Its entries $\tau^+_{ij}(u)\in\Yht[[u]]$  are found by
\beq\label{ltrn}
\tau^+_{ij}(u)=\delta_{ij}-h\sum_{s \geqslant 0} l_{ij}^+(-s)\ts e^{-su}
=\delta_{ij}-h\sum_{l\geqslant 0}\left(\sum_{s\geqslant 0} \frac{(-s)^l}{l!}l_{ij}^+(-s)\right)u^l.
\eeq
As with the matrix $T^+(u)$ given by \eqref{tplusu},   $\Tc^+(u)$  is  of the form $1+O(h)$.
Also, in contrast with $L^+(u)$, its constant term $\Tc^+(0)$ is no longer lower triangular; see \eqref{ltrn}. Finally,    \eqref{uqrtt3} and \eqref{equqrtt1} imply that $\Tc^+(u)$ satisfies the defining relations   for the  quantized universal enveloping algebra $\U $, so that we have
\begin{gather*}
R^+(e^{-u+v},e^{h/2})\ts \Tc_1^+(u)\ts \Tc_2^+(v) =\Tc_2^+(v)\ts \Tc_1^+(u)\ts R^+(e^{-u+v},e^{h/2}), \\
\Tc^{+}(u)\ts M \ts \Tc^{+}(u+ \kp h)^t\ts M^{-1} =1. 
\end{gather*} 
\end{rem}

\subsubsection{Annihilation operators}\label{subsec0302}
Let $\vac$ be the unit in the algebra $\U$. 
In the next lemma we construct the so-called annihilation operators for   types $B$, $C$, $D$; cf. \cite[Lemma 2.1]{EK5}.  

\begin{lem}\label{31}
For any $c\in\CC$ there exists a unique operator series 
$$
T^- (u)\in\ndo\CC^N \ot (\ndo\U )((u))_h
$$
such that $T^-(u)\vac=1\ot \vac$ and such that for all integers $k\geqslant 1$ we have
\begin{align}
&T^-_{k+1}(u)\ts T_1^+(v_1)\ldots T_k^+(v_k)\ts
=\, \R_{1\ts k+1}( -u+v_1-hc/2 ) \ldots \R_{ k\ts k+1}( -u+v_k-hc/2 )\non\\
&\qquad\times 
T_1^+(v_1)\ldots T_k^+(v_k)\R_{k\ts k+1}( -u+v_k+hc/2 )^{-1}\ldots \R_{1\ts k+1}( -u+v_1+hc/2 )^{-1}.\label{Th}
\end{align}
 Moreover, the   series $T^-(u)$ satisfies the identities
\begin{gather}
\R( u-v )\ts T_1^- (u)\ts T_2^- (v)=T_2^- (v)\ts T_1^- (u)\ts \R( u-v ),\label{RTThm}\\
T^- (u)\ts M\ts T^- (u+\kp h)^t\ts M^{-1} =1.\label{RELhm}
\end{gather} 
\end{lem}

\begin{prf}
First, we prove that the operator series $T^-(u)$ is well-defined by \eqref{Th}. As the coefficients of the matrix entries of all $T_1^+(v_1)\ldots T_k^+(v_k)$, along with $\vac$, span an $h$-adically dense $\CC[[h]]$-submodule of $\U$, it is sufficient to check that $T^-(u)$ preserves the ideal of its defining relations. For  relations \eqref{RTTh}, this is  verified by using the Yang--Baxter equation \eqref{YBEh} and arguing as in the proof of \cite[Lemma 2.1]{EK5}. As for   \eqref{RELh}, 
observe that
\begin{align}
T_2^-(u)\ts  \left(T_1^+(v) \ts  M_1\ts  T_1^+(v+\kp h)^t  -M_1\right)
=\R (-u+v-hc/2)\ts  T_1^+(v)\ts  X  -M_1,\label{rhs4}
\end{align}
where  
\begin{align*}
 X
=& \left(\R (-u+v-hc/2+ \kp h)\ts  T_1^+(v+h\kp)\right)^{t_1}\cdotrl Z\Fand\\
Z=&\left(\R (-u+v+hc/2+ \kp h)^{-1}\right)^{t_1}\cdotrl \left(\R (-u+v+hc/2)^{-1}M_1\right).
\end{align*}
By the crossing symmetry property \eqref{CSh} we have $Z=M_1$, so that
\begin{align*}
X=M_1 \ts   T_1^+(v+ \kp h)^t \ts \R (-u+v-hc/2+ \kp h)^{t_1}.
\end{align*}
Therefore, the right hand side of \eqref{rhs4} equals
$$
\R (-u+v-hc/2)\left(  T_1^+(v)\ts  M_1 \ts   T_1^+(v+ \kp h)^t \right) \R (-u+v-hc/2+ \kp h)^{t_1} -M_1. 
$$
However, by the crossing symmetry property \eqref{CSh} this is equal to
$$
\R (-u+v-hc/2)\left(  T_1^+(v)\ts  M_1 \ts   T_1^+(v+ \kp h)^t -M_1\right) \R (-u+v-hc/2+ \kp h)^{t_1} , 
$$
so it is clear that $T^-(u)$ maps \eqref{RELh} to the ideal of defining relations. This argument is easily generalized to an arbitrary element of 
the ideal of defining relations, thus implying that the operator series $T^-(u)$ is well-defined by \eqref{Th}. 
Moreover, as the $R$-matrix \eqref{Rh} belongs to $\ndo\CC^N\ot\ndo\CC^N ((u))[[h]]$, we conclude   that the  series $
T^-(u)$ is an element of $\ndo\CC^N \ot (\ndo\U )((u))_h
$, as required.

To complete the proof, it remains to verify the identities  \eqref{RTThm} and \eqref{RELhm}. First, we prove \eqref{RTThm}.
By applying  $\R_{k+1\ts k+2}(u_1-u_2)T^-_{k+1} (u_1)T^-_{k+2} (u_2)$  to 
$ T_{[k]}^+(v)\coloneqq T_1^+(v_1)\ldots T_k^+(v_k)$ with $  v=(v_1,\ldots ,v_k)$ and $k> 0$ an arbitrary integer,  we get
\begin{align*}
\R_{k+1\ts k+2}( u_1-u_2 ) \R^+_{1\ts k+2}\ldots \R^+_{k\ts k+2}\R^+_{1\ts k+1}\ldots \R^+_{k\ts k+1}
  T_{[k]}^+(v)
\R^{-1}_{k\ts k+1}\ldots \R^{-1}_{1\ts k+1} \R^{-1}_{k\ts n+2}\ldots \R^{-1}_{1\ts k+2},
\end{align*}
where
$$ 
\R_{i\ts k+j}^+=\R_{i\ts k+ j}( -u_j+v_i-hc/2 )
\fand  
\R_{i\ts k+j}^{-1}=\R_{i\ts k+ j}( -u_j+v_i+hc/2 )^{-1}.
$$
Using the Yang--Baxter equation \eqref{YBEh} we rewrite this expression as
\begin{align*}
\R^+_{1\ts k+1}\ldots \R^+_{k\ts k+1}
 \R^+_{1\ts k+2}\ldots \R^+_{k\ts k+2}
  T_{[k]}^+(v)
 \R^{-1}_{k\ts k+2}\ldots \R^{-1}_{1\ts k+2}
\R^{-1}_{k\ts k+1}\ldots \R^{-1}_{1\ts k+1}
\R_{k+1\ts k+2}( u_1-u_2 ).
\end{align*}
Finally, we observe that this coincides with the action of 
$T^-_{k+2} (u_2) T^-_{k+1} (u_1)\R_{k+1\ts k+2}( u_1-u_2 )$ on $T_{[k]}^+(v)$, 
so that the  relation \eqref{RTThm} follows.

Let us verify the remaining identity   \eqref{RELhm}. As before, we find 
\begin{align*}
&T^-_{k+1} (u_1)T^-_{k+2} (u_2) T_{[k]}^+(v)=  \R^+_{1\ts k+2}\ldots \R^+_{k\ts k+2}\\
 &\qquad\times \R^+_{1\ts k+1}\ldots \R^+_{k\ts k+1}
  T_{[k]}^+(v)
\R^{-1}_{k\ts k+1}\ldots \R^{-1}_{1\ts k+1} \R^{-1}_{k\ts k+2}\ldots \R^{-1}_{1\ts k+2}.
\end{align*}
Applying the transposition to the $(k+2)$-nd tensor factor  and then conjugating the equality by $M_{k+2}=1^{\ot (k+1)}\ot M$ we get
\begin{align*}
&T^-_{k+1} (u_1) M_{k+2} T^-_{k+2} (u_2)^{t_{k+2}} M_{k+2}^{-1} T_{[k]}^+(v)
=M_{k+2} \left( \R^{-1}_{k\ts k+2}  \ldots  \R^{-1}_{1\ts k+2} \right)^{t_{k+2}} \\
&\qquad\cdotrl \left( \left( \R^+_{1\ts k+2} \ldots  \R^+_{k\ts k+2}\right)^{t_{k+2}}\R^+_{1\ts k+1}\ldots \R^+_{k\ts k+1}
  T_{[k]}^+(v)
\R^{-1}_{k\ts k+1}\ldots \R^{-1}_{1\ts k+1} \right)M_{k+2}^{-1}.
\end{align*}
Finally, by applying the multiplication $a\ot b\mapsto ab$ on the tensor factors $k+1$ and $k+2$ and setting $u_2=u_1+\kp h$ we obtain
\begin{align}
&\left( T^-_{k+1} (u_1) M_{k+1} T^-_{k+1} (u_1+\kp h)^{t_{k+1}} M_{k+1}^{-1}\right) T_{[k]}^+(v)
= \left( \R^{'+}_{1\ts k+1} \ldots  \R^{'+}_{k\ts k+1}\right)^{t_{k+1}} \label{rhs}\\
&\qquad\cdotlr \left(\R^+_{1\ts k+1}\ldots \R^+_{k\ts k+1}
  T_{[k]}^+(v)
\R^{-1}_{k\ts k+1}\ldots \R^{-1}_{1\ts k+1} M_{k+1}  \left( \R^{'-1}_{k\ts k+1}  \ldots  \R^{'-1}_{1\ts k+1} \right)^{t_{k+1}} \right)M_{k+1}^{-1} ,\non
\end{align}
where $M_{k+1}^{\pm 1}=1^{\ot k }\ot M^{\pm 1}$,
$$
\R_{i\ts k+1}^{'+}=\R_{i\ts k+ 1}( -u_1+v_i-hc/2- \kp h)
\fand
\R_{i\ts k+1}^{'-1}=\R_{i\ts k+ 1}( -u_1+v_i+hc/2 - \kp h)^{-1}.
$$
However, crossing symmetry   \eqref{CSh}  and unitarity   \eqref{Uh}   imply the equality
\beq\label{csuni}
\R^{-1}_{k\ts k+1}\ldots \R^{-1}_{1\ts k+1} M_{k+1}  \left( \R^{'-1}_{k\ts k+1}  \ldots  \R^{'-1}_{1\ts k+1} \right)^{t_{k+1}} =M_{k+1}.
\eeq
Thus, we conclude that the right hand side of \eqref{rhs} is equal to
\begin{align}
&\left( \R^{'+}_{1\ts k+1} \ldots  \R^{'+}_{k\ts k+1}\right)^{t_{k+1}}\cdotlr \left(\R^+_{1\ts k+1}\ldots \R^+_{k\ts k+1}
  T_{[k]}^+(v) M_{k+1}
\right)M_{k+1}^{-1} \non\\
&\qquad=\left( \R^{'+}_{1\ts k+1} \ldots  \R^{'+}_{k\ts k+1}\right)^{t_{k+1}}\cdotlr \left(\R^+_{1\ts k+1}\ldots \R^+_{k\ts k+1} M_{k+1}\right)
 M_{k+1}^{-1} T_{[k]}^+(v) 
.\label{rhs2}
\end{align}
Note that the equality in \eqref{rhs2} holds because $T_{[k]}^+(v)$ commutes with $M_{k+1}^{\pm 1}$ and, also, with the $(k+1)$-th tensor factor of $ ( \R^{'+}_{1\ts k+1} \ldots  \R^{'+}_{k\ts k+1} )^{t_{k+1}}$.
As with \eqref{csuni}, the crossing symmetry   \eqref{CSh}  and unitarity   \eqref{Uh} imply the identity
$$
\left( \R^{'+}_{1\ts k+1} \ldots  \R^{'+}_{k\ts k+1}\right)^{t_{k+1}}\cdotlr \left(\R^+_{1\ts k+1}\ldots \R^+_{k\ts k+1} M_{k+1}\right)=M_{k+1},
$$
so we conclude by \eqref{rhs2} that the right hand side of \eqref{rhs} equals $T_{[k]}^+(v)$. Hence we have proved 
$$
\left( T^-_{k+1} (u_1) M_{k+1} T^-_{k+1} (u_1+\kp h)^{t_{k+1}} M_{k+1}^{-1}\right) T_{[k]}^+(v)=T_{[k]}^+(v).
$$
With integer $k$ being arbitrary, this implies \eqref{RELhm}, as required.
\end{prf}

Using  the properties of the $R$-matrix \eqref{Rh} one easily checks that   the operator series $T^- (u)$  is invertible. Furthermore,     its inverse is found by
\begin{align*}
&T^-_{k+1}(u)^{-1}T_1^+(v_1)\ldots T_k^+(v_k)\\
=& \,M_1\ldots M_k \left( \R_{k\ts k+1}( -u+v_k-hc/2- \kp h)^{t_k} \ldots \R_{ 1\ts k+1}( -u+v_1-hc/2- \kp h)^{t_1} \right)\non\\
& \cdotlr \left(
M_1^{-1} \ldots M_k^{-1}
T_1^+(v_1)\ldots T_k^+(v_k)\R_{1\ts k+1}( -u+v_1+hc/2 ) \ldots \R_{k\ts k+1}( -u+v_k+hc/2 ) \right) . 
\end{align*}

The matrix  $T^+(u)$ can be   regarded as an element of $\ndo\CC^N \ot (\ndo\U )[[u]]$, where its action is given by the algebra multiplication.  Formula \eqref{Th} implies the following relation among the operators $T^\pm (u)$:
\beq\label{RTThp}
\R ( -v+u-hc/2 )\ts T_1^+(u)\ts T_2^-(v)=T_2^-(v)\ts T_1^+(u)\ts \R ( -v+u+hc/2 ).
\eeq
Due to unitarity property \eqref{Uh}, this relation can be equivalently written as
\beq\label{RTThp2}
\R ( v-u+hc/2 )\ts T^-_1(v)\ts T_2^+(u)=T_2^+(u)\ts T^-_1(v)\ts \R ( v-u-hc/2 ).
\eeq

From now on we consider   the underlying $\CC[[h]]$-module structure of $\U$, along with the corresponding operators $T^\pm (u)$, and we denote  it by  $\V$. The subscript $c\in\CC$ indicates  the action of $T^-(u)$, as given by \eqref{Th}.
The $RTT$-formulas \eqref{RTTh}, \eqref{RTThm} and \eqref{RTThp} can be   generalized   as follows. For any  $a\in\CC$, integers $m,k\geqslant 1$ and the families of variables $u=(u_1,\ldots ,u_k)$, $v=(v_1,\ldots ,v_m)$ introduce the series  
\beq\label{arrows}
\R_{km}^{12}( u-v+ah )= \prod_{i=1,\dots,k}^{\longrightarrow} 
\prod_{j=k+1,\ldots,k+m}^{\longleftarrow} \R_{ij}( u_i -v_{j-k}+ah ) 
 \eeq
with the arrows indicating the order of the factors.
We shall use the superscripts ${}^{t_i}$ with $i=1,2$ to indicate that the first or   second tensor factors of the $R$-matrices in the given product are transposed. For example, we have
\begin{align*}
\R_{km}^{12,t_1}( u-v+ah )= \prod_{i=1,\dots,k}^{\longrightarrow} 
\prod_{j=k+1,\ldots,k+m}^{\longleftarrow} \R_{ij}^{t_i},\quad
\R_{km}^{12,t_2}( u-v+ah )= \prod_{i=1,\dots,k}^{\longrightarrow} 
\prod_{j=k+1,\ldots,k+m}^{\longleftarrow} \R_{ij}^{t_j}
,
 \end{align*}
where $\R_{ij}=\R_{ij}( u_i -v_{j-n}+ah )$.
Also, for $a\in \CC$ and $k\geqslant 2$ we use the   abbreviations
$$
T_{[k]}^\pm(u)=T_1^\pm(u_1)\ldots T_k^\pm(u_k)
,\quad
T^\pm_{[k]} (z+u+ah)=T^\pm_1 (z+u_1+ah)\ldots T^\pm_k (z+u_k+ah)
$$
 for operators on $(\ndo\CC^N)^{\ot k} \ot \V$.
As with the type $A$ case  \cite[Eq. (2.9)]{EK5}, using the original relations \eqref{RTTh}, \eqref{RTThm} and \eqref{RTThp}, one obtains
the identities
\begin{align}
\R_{km}^{12}( u-v )T_{[k]}^{\pm 13}(u)T^{\pm 23}_{[m]}(v)
=&\,\,T_{[m]}^{\pm 23}(v)T_{[k]}^{\pm 13}(u)\R_{km}^{12}( u-v ),\label{rrt1}\\
\R_{km}^{12}( -v+u-hc/2 )T_{[k]}^{+13}(u)T_{[m]}^{-23} (v)
=&\,\,T_{[m]}^{-23} (v)T_{[k]}^{+13}(u)\R_{km}^{12}( -v+u+hc/2 ),\label{rrt2}
\end{align}
which hold for the operators on the tensor product
\beq\label{supers}
\smalloverbrace{(\ndo\CC^N)^{\ot k}}^{1} \ot \smalloverbrace{(\ndo\CC^N)^{\ot m}}^{2}\ot  \smalloverbrace{V}^{3},
\eeq
where $V=\V$ and the superscripts indicate the tensor factors as above.

\begin{rem}\label{remark32}
In this remark we discuss relation between the definition  of the  operator series $T^-(u)$, as given by \eqref{Th}, and the quantum affine algebra structure.
Let $V_c(\g_N)$ be the quotient  of $U_h(\wht{\g}_N)_c$ 
by the $h$-adically complete  left ideal
generated by all elements $l^-_{ij}(r)$, where $i,j=1,\ldots ,N$ and $r=0,1,\ldots .$ By \eqref{guqrtt2}
there exists an operator series $L^-(y)$  on $V_c(\g_N)$ such that   for all $k\geqslant 1$ and the variables $x=(x_1,\ldots ,x_k)$ we have
\beq\label{formulla1}
 L^{-23 }(y)\ts L_{[k]}^{+13}(x) 
=
\RR_{k1}^{12}(xe^{hc/2}/y) \ts L_{[k]}^{+13}(x) \ts \RR_{k1}^{12}(xe^{-hc/2}/y)^{-1},
\eeq
where    the superscripts indicate the tensor factors   as in \eqref{supers} with $m=1$ and $V=V_c(\g_N)$ and we use the same notation for the elements of the quantum affine algebra and for their images in the quotient $V_c(\g_N)$.
 As in Remark \ref{remark31}, suppose that  $V_c(\g_N)$ is suitably completed so that we can apply the substitutions
  $x_i=e^{-v_i}$ with $i=1,\ldots ,k$ to \eqref{formulla1}, thus getting
\beq\label{formulla5}
 L^{-23 }(y)\ts \Tc^{+13}_{[k]}(v)
=
\RR_{k1}^{12}(e^{-v+hc/2}/y) \ts \Tc^{+13}_{[k]}(v) \ts \RR_{k1}^{12}(e^{-v-hc/2}/y)^{-1},  
\eeq
where $\Tc^+_{[k]}(v)=L^+_1(e^{-v_1})\ldots L^+_k(e^{-v_k})$ and the $R$-matrix products are defined as in \eqref{formulla3},  
$$
\RR_{k1}^{12}(  e^{-v\pm hc/2}/y )=   \RR_{1\ts k+1}(   e^{-v_1\pm hc/2}/y )\ldots \RR_{k\ts k+1}(   e^{-v_k\pm hc/2}/y ).
$$
The $R$-matrices on the right hand side belong to 
$(\ndo\CC^N)^{\ot 2}[[v_1,\ldots, v_k,h,y^{-1}]]$;
recall \eqref{R}, \eqref{uk6} and \eqref{rrr1}.
Furthermore, by Proposition \ref{imp_prop}, for any choice of integers $a_1,\ldots ,a_k ,l>0$ there exists an integer $r\geqslant 0$
such that the coefficients of all monomials
\beq\label{formulla4}
v_1^{a'_1}\ldots v_k^{a'_k}   h^{l'},\qquad \text{where} \qquad 0\leqslant a'_{i} < a_i  \fand l'< l,
\eeq
in 
$(1-y)^r\RR_{k1}^{12}(  e^{-v\pm hc/2}/y )^{\pm 1}$ belong to $(\ndo\CC^N)^{\ot 2}[y^{\pm 1}]$ and such that the coefficients of all monomials \eqref{formulla4} in
$$
\left((1-y)^r\RR_{k1}^{12}(  e^{-v\pm hc/2}/y )^{\pm 1}\right)\Big|_{y=e^{-u}}^{\text{mod }v_1^{a_1},\ldots ,v_k^{a_k} ,  h^{l} }
\Fand
(1-e^{-u})^{r}\ts
\R_{k1}^{12}(-u+v\mp hc/2)^{\pm 1} 
$$
coincide.
  Observe that
we can not  simply apply the substitution $y=e^{-u}$ to the right hand side of \eqref{formulla5}  as the resulting expression does not need to be  defined.
 However, by the preceding discussion,   the substitution $y=e^{-u}$ in  
$$
(1-e^{-u})^{-2r}
\left((1-y)^{2r}\ts \RR_{k1}^{12}(e^{-v+hc/2}/y) \ts \Tc^{+13}_{[k]}(v) \ts \RR_{k1}^{12}(e^{-v-hc/2}/y)^{-1}\right)\Big|_{y=e^{-u}}^{\text{mod }v_1^{a_1},\ldots ,v_k^{a_k} ,  h^{l}},\Big.
$$
where $(1-e^{-u})^{-2r}$ is regarded as an element of $\CC((u))$,
is well-defined.
Finally, we observe that  
   the definition \eqref{Th}  of the  operator series $T^-(u)$ admits the same form. Indeed,   by the choice of the integer $r$,  the coefficients of all monomials
 \eqref{formulla4} in
$$
(1-e^{-u})^{-2r}
\left((1-y)^{2r}\ts \RR_{k1}^{12}(e^{-v+hc/2}/y) \ts T^{+13}_{[k]}(v) \ts \RR_{k1}^{12}(e^{-v-hc/2}/y)^{-1}\right)\Big|_{y=e^{-u}}^{\text{mod }v_1^{a_1},\ldots ,v_k^{a_k} ,  h^{l}},\Big.
$$
and in \eqref{Th} coincide.
Also, it is worth noting that such form of the definition is in tune with the   products of $h$-adically quasi-compatible vertex operators; see  \cite[Sect. 4]{Li}. 
\end{rem}

\subsubsection{Braiding}\label{subsec0303}

From now on, the tensor products of $\CC[[h]]$-modules are understood as $h$-adically completed. 
Consider the $\CC[[h]]$-module map $\Dc\colon\V\to\V$  defined by 
\beq\label{D}
\Dc \vac =0\Fand\Dc\ts  T_{[k]}^+(u_1,\ldots ,u_k)\vac =
\sum_{i=1}^k \frac{\partial}{\partial u_i} T_{[k]}^+(u_1,\ldots ,u_k)\vac \quad\text{for}\quad k\geqslant 1.
\eeq
In the following lemma, we construct the braiding map. We use the arrow  in the superscript to indicate the opposite order  of the factors, for example we write (cf. \eqref{arrows})
$$
\R_{mk}^{\overleftarrow{12},t_1}(z+ u-v+ah )= \prod_{i=1,\dots,m}^{\longleftarrow} 
\prod_{j=m+1,\ldots,m+k}^{\longrightarrow} \R_{ij}( z+u_i -v_{j-k}+ah )^{t_i} .
$$

\begin{lem}\label{32}
There exists a unique  $\CC[[h]]$-module map
\beq\label{Smap2}
\Sc (z)\colon \V\ot\V\to\V\ot\V\ot\CC((z))[[h]]
\eeq
such that for all $m,k\geqslant 0$ and the variables
$u=(u_1,\ldots ,u_m) $ and $v=(v_1,\ldots ,v_k)$ we have
\begin{align}
&\Sc^{34}(z)\left(T_{[m]}^{+13}(u)T_{[k]}^{+24}(v)(\vac\ot\vac)\right)
=
\left(
M_{[m]}^1
\R_{mk}^{\overleftarrow{12},t_1}( z+u-v-h(c+\kp) )  
(M_{[m]}^1)^{-1}
\right) \label{Smap}\\
&\qquad \cdotlr
\Big(\R_{mk}^{12}( z+u-v )
   T_{[m]}^{+13}(u)
\R_{mk}^{12}( z+u-v+hc )^{-1}T_{[k]}^{+24}(v)\R_{mk}^{12}( z+u-v )(\vac\ot\vac)\Big),\non
\end{align}
where $M_{[m]}^1=M^{\ot m} \ot 1^{\ot k}$ and the superscripts   indicate the tensor factors as follows:
\beq\label{tfactors}
\smalloverbrace{(\ndo\CC^N)^{\ot m}}^{1}
\ot
\smalloverbrace{(\ndo\CC^N)^{\ot k}}^{2}
\ot 
\smalloverbrace{\V}^{3}
\ot 
\smalloverbrace{\V}^{4}.
\eeq
Moreover, the map $\Sc(z)$ satisfies 
the shift condition
\beq\label{shift}
\left[\Dc\ot 1,\Sc(z)\right]=-\frac{d}{dz}\Sc(z),
\eeq 
the Yang--Baxter equation
\beq\label{YBES}
\Sc_{12}(z_1)\ts \Sc_{13}(z_1+z_2)\ts \Sc_{23}(z_2)
=\Sc_{23}(z_2)\ts \Sc_{13}(z_1+z_2)\ts \Sc_{12}(z_1)
\eeq
and   the unitarity condition
\beq\label{US}
\Sc_{12}(z)\ts \Sc_{21}(-z)=1.
\eeq
\end{lem}

\begin{prf}
The fact that \eqref{Smap} defines a unique $\CC[[h]]$-module map \eqref{Smap2} is verified by an argument which goes in parallel with the corresponding part of the proof of \cite[Thm. 2.2]{BJK}. It relies on the properties of the $R$-matrix \eqref{Rh}, as given by Proposition \ref{Rpro}. In particular, note that the coefficients of matrix entries of all
$T_{[m]}^{+13}(u)T_{[k]}^{+24}(v)(\vac\ot\vac)$ are $h$-adically dense in $\V\ot\V$, which implies uniqueness. The shift condition \eqref{shift} is also proved by following the argument from \cite[Thm. 2.2]{BJK}.
Due to the crossing symmetry property \eqref{CSh}, formula \eqref{Smap} can be equivalently written as
\begin{align}
&\Sc^{34}(z) \ts  \R_{mk}^{12}( z+u-v )^{-1}\ts T_{[k]}^{+24}(v)\ts \R_{mk}^{12}( z+u-v-hc )\ts T_{[m]}^{+13}(u)(\vac\ot\vac)\non 
 \\
&\qquad \qquad=
   T_{[m]}^{+13}(u)\ts 
\R_{mk}^{12}( z+u-v+hc )^{-1}\ts T_{[k]}^{+24}(v)\ts \R_{mk}^{12}( z+u-v )(\vac\ot\vac). \label{pom2}
\end{align}
Finally, the proof of the Yang--Baxter equation \eqref{YBES} and the unitarity condition \eqref{US} is carried out 
by employing the corresponding properties of the $R$-matrix, \eqref{YBEh} and \eqref{Uh}, along with \eqref{pom2}, and arguing
as in the proof \cite[Thm. 4.1]{JKMY}. 
\end{prf}

\begin{rem}
As with the type $A$ case   \cite[Sect. 2]{EK5},  the braiding   is defined so that the $\Sc$-locality property \eqref{sloc} holds.  Its form  \eqref{Smap}, see also \eqref{pom2}, comes from the additive version of the so-called quantum current commutation relation which goes back to  \cite{RS}.  This relation, which holds for the  operator series $T(u)\coloneqq T^+(u)T^-(u+hc/2)^{-1}$  defining the vertex operator map \eqref{Y}, takes the form
\begin{align*}
&T_1 (u) \ts \R_{12}(u-v+hc)^{-1}\ts T_2 (v) \ts \R_{12}(u-v) =\R_{12}(-v+u)^{-1}\ts T_2 (v) \ts \R_{12}(-v+u-hc)\ts  T_1  (u) .
\end{align*}
\end{rem}

\subsubsection{Vertex operator map}\label{subsec0304}

Finally, we construct the vertex operator map $Y(\cdot ,z)$ which, along with the braiding $\Sc(z)$, defines a structure of $h$-adic quantum vertex algebra   over $\V$. The notion of  $h$-adic quantum vertex algebra was introduced by Li \cite[Def. 2.20]{Li} by generalizing the notion of quantum VOA of Etingof and Kazhdan \cite[Sect. 1.4.1]{EK5}.

\begin{thm}\label{thm36}
There exists a unique  $\CC[[h]]$-module map
\begin{align*}
Y\colon\V\ot\V&\to\V((z))_h\\
u\ot v&\mapsto Y(z)(u\ot v) = Y(u,z)v=\sum_{r\in \ZZ} u_r v\ts z^{-r-1}
\end{align*}
such that for all $k\geqslant 0$ and the variables $u=(u_1,\ldots ,u_k)$ we have
\beq\label{Y}
Y(T_{[k]}^+(u )\vac, z)
= T_{[k]}^+(z+u )
T_{[k]}^-(z+u+hc/2 )^{-1}.
\eeq
In particular, $Y(\vac,z)=1$. Moreover,  the given vertex operator map
satisfies the   weak associativity:
for any $u,v,w\in \V$ and $k\in\mathbb{Z}_{\geqslant 0}$
there exists $r\in\mathbb{Z}_{\geqslant 0}$
such that
\begin{align}
&(z_0 +z_2)^r\ts Y(u,z_0 +z_2)Y(v,z_2)\ts w \non \\
&\qquad  - (z_0 +z_2)^r\ts Y\big(Y(u,z_0)v,z_2\big)\ts w
\in h^k \V[[z_0^{\pm 1},z_2^{\pm 1}]], \label{wassoc}
\end{align}
   the $\mathcal{S}$-locality:
for any $u,v\in \V$ and $k\in\mathbb{Z}_{\geqslant 0}$ there exists
$r\in\mathbb{Z}_{\geqslant 0}$ such that   for all $w\in \V$
\begin{align}
&(z_1-z_2)^{r}\ts Y(z_1)\big(1\otimes Y(z_2)\big)\big(\mathcal{S}(z_1 -z_2)(u\otimes v)\otimes w\big) \non
 \\
 &\qquad-(z_1-z_2)^{r}\ts Y(z_2)\big(1\otimes Y(z_1)\big)(v\otimes u\otimes w)
 \in h^k \V[[z_1^{\pm 1},z_2^{\pm 1}]]\label{sloc}
\end{align}
and the   hexagon identity  
\begin{align} 
\Sc(z_1)\left(Y(z_2)\ot 1\right) =\left(Y(z_2)\ot 1\right)\Sc_{23}(z_1)\Sc_{13}(z_1+z_2).\label{hex}
\end{align}
Hence $(\V,Y,\vac, \Sc)$ is an $h$-adic quantum vertex algebra.
\end{thm}

\begin{prf}
The proof that the vertex operator map is well-defined by \eqref{Y}
relies on the defining relations
\eqref{RTTh} and \eqref{RELh} for the quantized universal enveloping algebra and also on \eqref{RTThm} and \eqref{RELhm}.
It goes  similarly to the first part of the proof of Lemma \ref{31}. The weak associativity \eqref{wassoc},    $\Sc$-locality \eqref{sloc}  and the hexagon identity \eqref{hex}  are verified by calculations that go in parallel with the rational $R$-matrix case and rely   on the properties of the trigonometric $R$-matrix established by Proposition \ref{Rpro}; see the proofs of  \cite[Thm. 2.2]{BJK}, \cite[Thm. 2.3.8]{G} and \cite[Thm. 4.1]{JKMY}. 
Moreover, we note that the map \eqref{D} satisfies $\Dc v=v_{-2}\vac$ for all $v\in\V$.
As for the last assertion of the theorem,   due to Lemmas \ref{31},  \ref{32} and the preceding discussion, $(\V,Y,\vac, \Sc)$   satisfies all   axioms   of   $h$-adic quantum vertex algebra, as given in \cite[Def. 2.20]{Li}.  
\end{prf}

\subsection{\texorpdfstring{$\phi$}{phi}-coordinated modules
}\label{subsec0305}
In this subsection, we recall the notion of $\phi$-coordinated $\V$-module and we show that any restricted $U_h(\wht{\g}_N)_c$-module 
is naturally equipped with the structure of $\phi$-coordinated $\V$-module with $\phi(z_2,z_0) = z_2 e^{-z_0}$.

\subsubsection{Main theorem
}\label{subsec03050305}
 The notion of $\phi$-coordinated   module, 
where $\phi$ is an associate of the one-dimensional additive formal group, was introduced by Li  \cite[Def. 3.4]{Li1}.
The following definition of $\phi$-coordinated module   is obtained as a straightforward
generalization of the original definition over the ring $\CC[[h]]$, which makes it  compatible with $h$-adic quantum vertex algebras.
Even though the definition is formulated in terms of   the particular  associate 
$ \phi(z_2,z_0) = z_2 e^{-z_0}$, 
  it can be easily generalized to  an arbitrary choice of an associate.

\begin{defn}\label{phimod}
Let $V$ be an $h$-adic quantum vertex algebra.
A {\em $\phi$-coordinated $V$-module}, where $ \phi(z_2,z_0) = z_2 e^{-z_0}$, is a pair $(W,Y_W)$ such that $W$ is a topologically free $\CC[[h]]$-module and $Y_W=Y_W(\cdot, z)$ is  a $\mathbb{C}[[h]]$-module map
\begin{align*}
Y_W \colon  V \ot W&\to W((z))_h\\
u\ot w&\mapsto Y_W(z)(u\ot w)=Y_W(u,z)w=\sum_{r\in\mathbb{Z}} u_r w \ts z^{-r-1} 
\end{align*}
which satisfies
 $Y_W(\vac,z)w=w$  for all $w\in W$ and  
  the {\em weak associativity}: for any $u,v\in  V $, $k\in\mathbb{Z}_{\geqslant 0}$ there exists $r\in\mathbb{Z}_{\geqslant 0}$ such that
\begin{align}
&(z_1-z_2)^r\ts Y_W(u,z_1)Y_W(v,z_2)\in\om (W,W((z_1,z_2)) )\mod h^k, \label{associativitymod0}\\
&\big((z_1-z_2)^r\ts Y_W(u,z_1)Y_W(v,z_2)\big)\big|_{z_1= z_2 e^{-z_0}}^{\modd h^k}  \big. \non\\
&\qquad- z_2^r(e^{-z_0} - 1)^r\ts Y_W\left(Y(u,z_0)v,z_2\right)\ts
\in\ts  h^k \om(W,W[[z_0^{\pm 1},z_2^{\pm 1}]]).\label{associativitymod}
\end{align}
Let $W_1$ be a topologically free $\CC[[h]]$-submodule of $W$. A pair $(W_1,Y_{W_1})$ is said to be a {\em $\phi$-coordinated $V$-submodule} of $W$ if $Y_W(v,z)w$ belongs to $W_1[[z^{\pm 1}]]$ for all $v\in V$ and $w\in W_1$, where $Y_{W_1}$ denotes the restriction and corestriction of $Y_W$,
$$Y_{W_1} (z)=Y_{W} (z)\big|_{\V \ot W_1}^{ W_1((z))_h}\Big.\colon V\ot W_1 \to W_1 ((z))_h.$$
\end{defn}

Regarding the definition, observe that the formulation  of the weak associativity axiom, \eqref{associativitymod0} and \eqref{associativitymod} employs the notation from the last paragraph of Section \ref{sec02}.

The following theorem is the main result of this paper. The substitutions $x_i=ze^{-u_i}$  on the right hand side of \eqref{Ywmap} are carried out simultaneously for all $i=1,\ldots ,k$.

\begin{thm}\label{main1}
Let $c\in \CC$ and let $W$ be a restricted $U_h(\wht{\g}_N)_c$-module.
There exists a unique structure of $\phi$-coordinated $\V$-module on $W$ such that for all $k\geqslant 1$ we have
\beq\label{Ywmap}
Y_W(T_{[k]}^+(u_1,\ldots ,u_k),z)=
L_{[k]}^+(x_1,\ldots, x_k)\big|_{x_i=ze^{-u_i}}
L_{[k]}^-(x_1 e^{-hc/2},\ldots, x_k e^{-hc/2})^{-1}\big|_{x_i=ze^{-u_i}}.
\eeq
Furthermore, if \eqref{Ywmap} defines a structure of irreducible $\phi$-coordinated $\V$-module on $W$, then $W$ is irreducible as an  $U_h(\wht{\g}_N)_c$-module as well.
\end{thm}

\begin{prf}
The first assertion of the theorem follows from Lemmas \ref{lmaa} and \ref{lmab}  given in Subsection \ref{subsec03050306} below. As for the second assertion, suppose $W_1\subseteq W$ is a restricted $U_h(\wht{\g}_N)_c$-submodule.    The coefficients of matrix entries of all $T_{[k]}^+(u_1,\ldots ,u_k)$ with $k\geqslant 1$ along with $\vac$ span an $h$-adically dense $\CC[[h]]$-submodule of $\V$, so  that  \eqref{Ywmap} implies     $Y(v,z)W_1\subseteq W_1[[z^{\pm 1}]]$ for all $v\in\V$. Thus,  $W_1$ is a $\V$-submodule of $W$.
\end{prf}

\begin{rem}\label{yremark}
The rational $R$-matrix counterpart of Theorem \ref{main1} can be established analogously. More specifically, consider the ($h$-adic) quantum vertex algebra $\V^{\text{rat}}$ associated with the rational $R$-matrix  in type $A$ \cite[Thm. 2.3]{EK5} or in types $B$, $C$, $D$ \cite[Thm. 2.2]{BJK}. In parallel with Section \ref{sec05}, one can introduce the notion of restricted module for the level $c\in\CC$ double Yangian in types $A$--$D$  defined over  $\CC[[h]]$, using its $R$-matrix presentation  \cite{I,JLY}. Then, by arguing as in the proof of Theorem \ref{main1}, one shows that such restricted modules are naturally equipped with the structure of $\V^{\text{rat}}$-module; see also \cite[Thm. 6.11]{Li} for the  double Yangian in type $A_1$. In particular, the rational $R$-matrix setting does not    require  the use of theory of $\phi$-coordinated modules, as the usual notion of $h$-adic quantum vertex algebra   module  \cite[Def. 2.23]{Li} is suitable. 
\end{rem}

\begin{rem}
In this remark, we discuss an approach in the context of quantum current algebras  which might lead to the proof of the converse of Theorem \ref{main1}.
In \cite{Koz1} we established an equivalence of certain  $\phi$-coordinated modules   for the Etingof--Kazhdan quantum affine vertex algebra in type $A$ \cite{EK5}  and restricted modules for (slightly modified) Ding's quantum current algebra   \cite{D} defined over $\CC[[h]]$. 
On the other hand, by using the famous Ding--Frenkel isomorphism \cite{DF}, Ding proved that  the original quantum current algebra over $\CC(q)$ is isomorphic to the quantum affine algebra in type $A$   \cite[Prop. 3.1]{D}. Due to the lack of such quantum current realizations for   quantum affine algebras in types $B$, $C$ and $D$, we were not able to prove the equivalence of $\phi$-coordinated modules for the   $h$-adic quantum  vertex algebras  established in Theorem \ref{thm36} and restricted modules for the quantum affine algebras of the same type. Instead, we only obtained a partial result, as given in Theorem \ref{main1}. Nonetheless, our calculations   suggest that it might be plausible to introduce quantum current algebras in types $B$, $C$ and $D$ via following defining relations. As with the type $A$, quantum current algebras would be defined as topological algebras over $\CC[[h]]$ in  generators  $\lambda_{ij}^{(r)}$, where $i,j=1,\ldots ,N$ and $r\in\mathbb{Z}$.  
One family of defining relations is given by the so-called quantum current commutation relation, cf. \cite{RS},
\beq\label{qcrem}
\Lc_1(x)\ts \RR_{21}(ye^{hc}/x)\ts\Lc_2(y)\ts \RR_{21}(y/x)^{-1}
=
\RR_{12}(x/y)^{-1}\ts \Lc_2(y)\ts \RR_{12}(xe^{hc}/y)\ts \Lc_1(x),
\eeq
where
$$
\Lc (x)=\sum_{i,j=1}^N e_{ij}\ot \lambda_{ij} (x) \quad\text{with}\quad
\lambda_{ij} (x)=\delta_{ij}- h\sum_{r\in\ZZ}  \lambda_{ij}^{(  r)} x^{- r-1}. 
$$
Denote by $\Lc_{[2]}(x,y)$  the left hand side of \eqref{qcrem}. Let  $m_{12}\colon a\ot b\mapsto ab$, $m_{21}\colon a\ot b\mapsto ba$ be the multiplications $(\ndo\CC^N)^{\ot 2}\to\ndo\CC^N$.
The other family of defining relations takes the form
$$
m_{ab}\left(M_b\ts \Lc_{[2]}(x,y)^{t_b}\big|_{y=x\zeta}\big.\ts M_b^{-1}\right)=1
\fand
m_{ab}\left(M_a\ts \Lc_{[2]}(x ,y)^{t_a}\big|_{y=x\zeta^{-1}}\big.\ts M_a^{-1}\right)=1 
$$ 
for $(a,b)=(1,2),(2,1)$, where
 the transposition is defined as in Section \ref{sec02}, the diagonal matrix $M $ is given by \eqref{CSh} and $\zeta = \xi\left|_{q=e^{h/2}}\right.$; recall \eqref{xi}.
Of course, it remains to see whether such quantum current   algebras would lead to new realizations of quantum affine algebras in types $B$, $C$ and $D$ and, furthermore, whether they would govern the $\phi$-coordinated representation theory of the   $h$-adic quantum  vertex algebras  from Theorem \ref{thm36}, as is the case in   type $A$.
\end{rem}

\subsubsection{Two lemmas}
 \label{subsec03050306}
This subsection is dedicated to the proof of Lemmas \ref{lmaa} and \ref{lmab}  which finalize  the proof of Theorem \ref{main1}.

\begin{lem}\label{lmaa}
Let $W$ be a restricted $U_h(\wht{\g}_N)_c$-module.
The formula  \eqref{Ywmap}, along with $Y_W(\vac,z)=1_W$, uniquely determines a $\CC[[h]]$-module map
$
  \V \to \om(W,W((z))_h).
$
\end{lem}

\begin{prf}
The coefficients of matrix entries of all $T_{[k]}^+(u_1,\ldots ,u_k)$ with $k\geqslant 0$ span a $h$-adically dense $\CC[[h]]$-submodule of $\V$. Hence, to prove the lemma, it is sufficient to check that $v\mapsto Y_W(v,z)$ preserves the ideal of defining relations \eqref{RTTh} and \eqref{RELh} of $\V$. 
First, we consider the relations \eqref{RTTh} for operators on the restricted $U_h(\wht{\g}_N)_c$-module $W$. By \eqref{equqrtt1} for any $k\geqslant 2$ and $i=1,\ldots ,k-1$ we  have
\begin{align*}
&R^+_{i\ts i+1}(x_i /x_{i+1} ,e^{h/2})\ts L_{[k]}^+(x )\ts 
L_{[k]}^-(  e^{-hc/2}  x )^{-1} \\
&\qquad -P_{i\ts i+1} L_{[k]}^+(x^{(i)})\ts 
L_{[k]}^-(  e^{-hc/2} x^{(i)})^{-1}\ts
P_{i\ts i+1}\ts
R^+_{i\ts i+1}(x_i /x_{i+1} ,e^{h/2})=0,
\end{align*}
where  
 $x=(x_1,\ldots ,x_k)$  and $x^{(i)}=(x_1,\ldots ,x_{i-1},x_{i+1},x_{i},x_{i+2},\ldots ,x_k).$ 
Applying the substitutions $x_i=ze^{-u_i}$ with $i=1,\ldots ,k$, the equality takes the form
\begin{align*}
&R^+_{i\ts i+1}(e^{-u_i+u_{i+1}} ,e^{h/2})\ts L_{[k]}^+(x )\big|_{x_i=ze^{-u_i}}\ts 
L_{[k]}^-(  e^{-hc/2}  x )^{-1} \big|_{x_i=ze^{-u_i}}\\
&\qquad -P_{i\ts i+1} L_{[k]}^+(x^{(i)})\big|_{x_i=ze^{-u_i}}\ts 
L_{[k]}^-(  e^{-hc/2} x^{(i)})^{-1}\big|_{x_i=ze^{-u_i}}\ts
P_{i\ts i+1}\ts
R^+_{i\ts i+1}(e^{-u_i+u_{i+1}} ,e^{h/2})=0.
\end{align*}
The left hand side of the equality coincides with the image of the left hand side of 
\begin{align*}
&R^+_{i\ts i+1}(e^{-u_i+u_{i+1}} ,e^{h/2})\ts T_{[k]}^+(u_1,\ldots ,u_k ) \\
&\qquad -P_{i\ts i+1} T_{[k]}^+(u_1,\ldots ,u_{i-1},u_{i+1},u_{i},u_{i+2},\ldots ,u_k) \ts
P_{i\ts i+1}\ts
R^+_{i\ts i+1}(e^{-u_i+u_{i+1}} ,e^{h/2})=0
\end{align*}
under $Y_W(\cdot,z)$. Thus, we conclude that \eqref{Ywmap} preserves the defining  relations \eqref{RTTh}, as they are equivalent to \eqref{eqRTTh}.

As for the remaining relation  \eqref{RELh}, to simplify notation, we  write $T_j^+ =T_j^+(u_j)$ and $\bar{T}_j^+ =T_j^+(u_j+h\kp)^t$. It is sufficient to check that  the image of
\beq\label{zero1}
T_1^+\ldots T_{i-1}^+ \ts T_i^+\ts M_i\ts  \bar{T}_i^+ \ts M_i^{-1}\ts T_{i+1}^+\ldots T_{n}^+ - T_1^+\ldots T_{i-1}^+ \ts   T_{i+1}^+\ldots T_{k}^+
\eeq
under $Y_W(\cdot,z)$ is zero for all $n\geqslant 1$ and $i=1,\ldots ,n$. Write 
$$L_j^{+ } = L_j^{+ }(x_j),\quad L_j^{- } = L_j^{- }(x_j e^{-hc/2})^{-1},\quad \bar{L}_j^{+ } = L_j^{+ }(x_j\zeta)^t,\quad \bar{L}_j^{- } = \left(L_j^{- }(x_j e^{-hc/2}\zeta)^{-1}\right)^{t}.$$  
 By applying $Y_W(\cdot,z)$ to \eqref{zero1} we get
\begin{align}
\left(
L_1^+\ldots   L_i^+\left( 
\bar{L}_i^+\cdotrl\left(
L_{i+1}^+\ldots L_k^+\ts 
 X\right) \right)M_i^{-1}
\right)\Big|_{x_i=ze^{-u_i}}-Z,\label{zero2}
\end{align}
where
\begin{align*}
&X=
\left( L_k^-  \ldots  L_{i+1}^-  
\left(
 \bar{L}_i^{-}  \cdotrl
\left( L_i^-  M_i\right)
\right)
 L_{i-1}^-  \ldots  L_1^- 
\right)\Big|_{x_i=ze^{-u_i}}, \\
&Z=
\left(
L_1^+\ldots   L_{i-1}^+   \ts   
L_{i+1}^+\ldots L_k^+\right)\Big|_{x_i=ze^{-u_i}}
 \left( L_k^-  \ldots  L_{i+1}^-  
 \ts   
 L_{i-1}^-  \ldots  L_1^- 
\right)\Big|_{x_i=ze^{-u_i}}.
\end{align*}
However, defining relation \eqref{uqrtt3} implies the identity
$$
X=
\left( L_k^-  \ldots  L_{i+1}^-  
 \ts   M_i\ts 
 L_{i-1}^-  \ldots  L_1^- 
\right)\Big|_{x_i=ze^{-u_i}}
= M_i\ts 
\left( L_k^-  \ldots  L_{i+1}^-  
 \ts   
 L_{i-1}^-  \ldots  L_1^- 
\right)\Big|_{x_i=ze^{-u_i}}
, 
$$
so that \eqref{zero2} takes the form
\begin{align*}
&\left(
L_1^+\ldots  L_{i-1}^+   \ts   L_i^+ \ts M_i \ts   
\bar{L}_i^+ 
L_{i+1}^+\ldots L_k^+\right)\Big|_{x_i=ze^{-u_i}}
 \left( L_k^-  \ldots  L_{i+1}^-  
 \ts   
 L_{i-1}^-  \ldots  L_1^- 
\right)\Big|_{x_i=ze^{-u_i}}
   M_i^{-1} -Z\\ 
=&	\left(
L_1^+\ldots  L_{i-1}^+   \ts   L_i^+ \ts M_i \ts   
\bar{L}_i^+ \ts M_i^{-1}
L_{i+1}^+\ldots L_k^+\right)\Big|_{x_i=ze^{-u_i}}
 \left( L_k^-  \ldots  L_{i+1}^-  
 \ts   
 L_{i-1}^-  \ldots  L_1^- 
\right)\Big|_{x_i=ze^{-u_i}}
    -Z.
\end{align*}
Using  \eqref{uqrtt3} once again, we rewrite the first term of the expression above, thus getting  
\begin{align*}
\left(
L_1^+\ldots   L_{i-1}^+   \ts   
L_{i+1}^+\ldots L_k^+\right)\Big|_{x_i=ze^{-u_i}}
 \left( L_k^-  \ldots  L_{i+1}^-  
 \ts   
 L_{i-1}^-  \ldots  L_1^- 
\right)\Big|_{x_i=ze^{-u_i}}
-Z=0, 
\end{align*}
as required.
Therefore, we have proved that
\eqref{Ywmap} uniquely defines a $\CC[[h]]$-module map
$ 
  \V \to (\ndo W)[[z^{\pm 1}]]$.
Finally, with $W$ being restricted, we conclude from \eqref{restr5} that its image belongs to
$\om(W,W((z))_h)$, thus finishing the proof.
\end{prf}

\begin{lem}\label{lmab}
Let $W$ be a restricted $U_h(\wht{\g}_N)_c$-module.
The $\CC[[h]]$-module  map  \eqref{Ywmap} possesses the weak associativity properties \eqref{associativitymod0} and \eqref{associativitymod}. 
\end{lem}

\begin{prf}
We start by verifying \eqref{associativitymod0}. For any integers $k,m\geqslant 1$ we consider   the expression 
\beq\label{spr3}
T_{[k]}^{+23}(u )\ts T_{[m]}^{+14}(v )=T_{[k]}^{+23}(u_1,\ldots,u_k)\ts  T_{[m]}^{+14}(v_1,\ldots ,v_m),
\eeq
whose coefficients belong to \eqref{tfactors}, 
where the superscripts $1,2,3,4$  indicate the tensor factors in accordance with \eqref{tfactors}. 
By applying $Y_W(z_1)(1\ot Y_W(z_2))$ on the given expression and using \eqref{Ywmap} we get
\begin{align}
&\, Y_W(T_{[n]}^{+23}(u) ,z_1 )Y_W(T_{[m]}^{+13}(v) ,z_2 )\label{exprr01}
\\
= &\,
L_{[k]}^{+23} (x) \big|_{x_i =z_1 e^{-u_i}}\big.  
L_{[k]}^{-23} (e^{-hc/2}x)^{-1} \big|_{x_i =z_1 e^{-u_i}}\big.  
L_{[m]}^{+13} (y) \big|_{y_j =z_2 e^{-v_j}}\big.  
L_{[m]}^{-13} (e^{-hc/2}y)^{-1} \big|_{y_j =z_2 e^{-v_j}}\big.   ,
\non
\end{align}
where $x=(x_1,\ldots ,x_k)$ and $y=(y_1,\ldots ,y_m)$. Let $(M_{[m]}^1)^{\pm 1}=M_1^{\pm 1}\ldots M_m^{\pm 1}$ and
$$
A_{mk}^{12}( z_2/z_1)= M_{[m]}^1 \prod_{j=1,\dots,m}^{\longleftarrow} 
\prod_{i=m+1,\ldots,m+k}^{\longrightarrow} \RR_{ji}( y_{j} e^{hc}\zeta^{-1}/x_{i-k})^{t_j}\bigg|_{\substack{\hspace{-18pt} y_j =z_2 e^{-v_j}\\x_{i-k} =z_1 e^{-u_{i-k}}}}   (M_{[m]}^1)^{- 1}.\bigg.
 $$
By the crossing symmetry relation \eqref{CSv2} we have
$$
A_{mk}^{12}( z_2/z_1 ) \cdotlr \RR_{mk}^{12}( ye^{hc}/x )\Big|_{\substack{y_j =z_2 e^{-v_j}\\x_i =z_1 e^{-u_i}}}=1.\Big.
$$
Hence, by using the $RLL$-relation \eqref{guqrtt2} we rewrite the right hand side of \eqref{exprr01} as
\begin{align}
A_{mk}^{12}( z_2/z_1 ) \cdotlr\bigg(
&L_{[k]}^{+23} (x) \big|_{x_i =z_1 e^{-u_i}}\big.   
L_{[m]}^{+13} (y) \big|_{y_j =z_2 e^{-v_j}}\big. 
\RR_{mk}^{12}( y /x )\Big|_{\substack{y_j =z_2 e^{-v_j}\\x_i =z_1 e^{-u_i}}} \Big.\bigg.\non\\ 
\bigg. &\times L_{[k]}^{-23} (e^{-hc/2}x)^{-1} \big|_{x_i =z_1 e^{-u_i}}\big. 
L_{[m]}^{-13} (e^{-hc/2}y)^{-1} \big|_{y_j =z_2 e^{-v_j}}\big.\bigg).\non
\end{align}
Next,  we   reorder the last three factors using \eqref{guqrtt1}, thus getting
\begin{align}
A_{mk}^{12}( z_2/z_1 ) \cdotlr\bigg(
&L_{[k]}^{+23} (x) \big|_{x_i =z_1 e^{-u_i}}\big.   
L_{[m]}^{+13} (y) \big|_{y_j =z_2 e^{-v_j}}\big. 
L_{[m]}^{-13} (e^{-hc/2}y)^{-1} \big|_{y_j =z_2 e^{-v_j}}
\bigg.\non\\ 
\bigg. &\times L_{[k]}^{-23} (e^{-hc/2}x)^{-1} \big|_{x_i =z_1 e^{-u_i}}\RR_{mk}^{12}( y /x )\Big|_{\substack{y_j =z_2 e^{-v_j}\\x_i =z_1 e^{-u_i}}} \Big.
\big. 
\big.\bigg).\label{assocpr1}
\end{align}

Let $a_1,\ldots ,a_k,b_1,\ldots ,b_m ,l>0$ be   integers. Consider the coefficients of the monomials 
\beq\label{monomials}
u_1^{a'_1}\ldots u_k^{a'_k} v_1^{b'_1}\ldots v_{m}^{b'_m} h^{l'},\qquad \text{where} \qquad 0\leqslant a'_{i} < a_i,\quad 0\leqslant  b'_j < b_j \fand l'< l,
\eeq
in \eqref{assocpr1}. By the first assertion of Proposition \ref{imp_prop}, there exists a nonnegative integer $r $ such that the coefficients of all monomials \eqref{monomials} in
$$
(z_1 -z_2)^r A_{mk}^{12}( z_2/z_1 ) \Fand 
(z_1 -z_2)^r \RR_{mk}^{12}( y /x )\Big|_{\substack{y_j =z_2 e^{-v_j}\\x_i =z_1 e^{-u_i}}} \Big.
$$
belong to $(\ndo\CC^N)^{\ot (m+k)}[z_1^{\pm 1},z_2^{\pm 1}]$. Furthermore, as $W$ is a restricted module, by \eqref{restr5}, the coefficients of all monomials \eqref{monomials} in
$$
L_{[m]}^{-13} (e^{-hc/2}y)^{-1} \big|_{y_j =z_2 e^{-v_j}}
  L_{[k]}^{-23} (e^{-hc/2}x)^{-1}\big|_{x_i =z_1 e^{-u_i}}w\quad\text{with}\quad w\in W\big.\big.
$$
belong to $(\ndo\CC^N)^{\ot (m+k)}\ot W [z_1^{- 1},z_2^{-1}]$. Therefore, we conclude from \eqref{assocpr1} that   for any $w\in W$ the coefficients of all monomials \eqref{monomials} in 
\begin{align}
(z_1 - z_2)^{2r} A_{mk}^{12}( z_2/z_1 ) \cdotlr\bigg(
&L_{[k]}^{+23} (x) \big|_{x_i =z_1 e^{-u_i}}\big.   
L_{[m]}^{+13} (y) \big|_{y_j =z_2 e^{-v_j}}\big. 
L_{[m]}^{-13} (e^{-hc/2}y)^{-1} \big|_{y_j =z_2 e^{-v_j}}
\bigg.\non \\ 
\bigg. &\times L_{[k]}^{-23} (e^{-hc/2}x)^{-1} \big|_{x_i =z_1 e^{-u_i}} w \ts \RR_{mk}^{12}( y /x )\Big|_{\substack{y_j =z_2 e^{-v_j}\\x_i =z_1 e^{-u_i}}} \Big.
\big. 
\big.\bigg)\label{exprnon} 
\end{align}
belong to $(\ndo\CC^N)^{\ot (m+k)}\ot W ((z_1,z_2))$, so that weak associativity property \eqref{associativitymod0}   follows.

It remains to prove \eqref{associativitymod}. Without loss of generality we can assume, due to the  second assertion of Proposition \ref{imp_prop}, that the integer $r$ has been chosen so that the coefficients of all monomials \eqref{monomials} in
$$ 
\left( (z_1 -z_2)^r A_{mk}^{12}( z_2/z_1 )\right)\big|_{z_1 =z_2 e^{-z_0}}^{\text{mod } u_1^{a_1},\ldots, u_k^{a_k} ,v_1^{b_1},\ldots, v_{m}^{b_m}, h^{l}} 
\big.
\fand
z_2^r (e^{-z_0}-1)^r B_{mk}^{12}(-z_0-u+v),
$$
where
$$ 
B_{mk}^{12}(-z_0-u+v) =   M_{[m]}^1  \prod_{j=1,\dots,m}^{\longleftarrow} 
\prod_{i=m+1,\ldots,m+k}^{\longrightarrow} \R_{ji}( -z_0-u_{i-k}+v_{j}-h(c+\kp))^{t_j}   (M_{[m]}^1)^{- 1},\bigg.
 $$
coincide. In addition, we can assume that   the coefficients of all monomials \eqref{monomials} in
$$
\left((z_1 -z_2)^r\RR_{mk}^{12}( y /x )\Big|_{\substack{y_j =z_2 e^{-v_j}\\x_i =z_1 e^{-u_i}}} \Big.\right)\bigg|_{z_1 =z_2 e^{-z_0}}^{\text{mod } u_1^{a_1},\ldots, u_k^{a_k} ,v_1^{b_1},\ldots, v_{m}^{b_m}, h^{l}} 
\bigg.\fand
z_2^r (e^{-z_0}-1)^r
\R_{mk}^{12}(-z_0 -u+v)
$$
coincide as well. Hence, by regarding  \eqref{exprnon} modulo $u_1^{a_1},\ldots, u_k^{a_k} ,v_1^{b_1},\ldots, v_{m}^{b_m}, h^{l}$ and then taking the substitution $z_1 =z_2 e^{-z_0}$ we obtain
\begin{align}
&z_2^{2r} (e^{-z_0}-1)^{2r} B_{mk}^{12}(-z_0-u+v) \cdotlr\bigg(
L_{[k]}^{+23} (x) \big|_{x_i =z_2 e^{-z_0-u_i}}\big.   
L_{[m]}^{+13} (y) \big|_{y_j =z_2 e^{-v_j}}\big. 
\bigg.  \non\\ 
&\bigg. 
\quad\times 
L_{[m]}^{-13} (e^{-hc/2}y)^{-1} \big|_{y_j =z_2 e^{-v_j}}
L_{[k]}^{-23} (e^{-hc/2}x)^{-1} \big|_{x_i =z_2 e^{-z_0-u_i}} w \ts \R_{mk}^{12}(-z_0 -u+v) \Big.
\big. 
\big.\bigg) \label{that2}
\end{align}
modulo $u_1^{a_1},\ldots, u_k^{a_k} ,v_1^{b_1},\ldots, v_{m}^{b_m}, h^{l}  $.
However, by a direct calculation which employs $RTT$-relations \eqref{rrt1} and \eqref{rrt2}, along with the definitions of the vertex operator map \eqref{Y} and   $\phi$-coordinated module map \eqref{Ywmap}, one   shows  that   \eqref{that2} is also obtained  by applying
$$
z_2^{2r} (e^{-z_0}-1)^{2r}\ts 
Y_W(z_2)(Y(z_0)\ot 1)
$$
to the expression
 $ T_{[k]}^{+23}(u )\ts T_{[m]}^{+14}(v )\ot w .$ 
Thus we conclude that the weak associativity property \eqref{associativitymod} holds.
\end{prf}

\section*{Acknowledgement}
This work has been supported in part by Croatian Science Foundation under the project UIP-2019-04-8488.

 \linespread{1.0}


\begin{thebibliography}{9}
\bibitem{BK}
B. Bakalov, V. G. Kac,
{\em Field algebras}, 
Int. Math. Res. Not. (2003), no. 3, 123--159;
\href{http://arxiv.org/abs/math/0204282}{arXiv:math/0204282 [math.QA]}.

\bibitem{BM}
C. Boyallian, V. Meinardi,
{\em An approach to Quantum Conformal Algebra}, 
\href{https://arxiv.org/abs/2012.15299v1}{arXiv:2012.15299 [math.QA]}.

\bibitem{BJK}
M. Butorac, N. Jing, S. Ko\v{z}i\'{c},
{\em $h$-Adic quantum vertex algebras associated with rational $R$-matrix in types $B$, $C$ and $D$},  Lett. Math. Phys. \textbf{109} (2019), 2439--2471;
\href{https://arxiv.org/abs/1904.03771}{arXiv:1904.03771 [math.QA]}.

\bibitem{CP}
 V. Chari, A. Pressley, 
{\em A guide to quantum groups}, Cambridge University Press,
Cambridge, 1994.

\bibitem{DGK}
A. De Sole, M. Gardini, V. G. Kac,
{\em On the structure of quantum vertex algebras},
J. Math. Phys. \textbf{61} (2020), 011701 (29pp);
\href{https://arxiv.org/abs/1906.05051}{arXiv:1906.05051 [math.QA]}.

\bibitem{D}
J. Ding, 
{\em Spinor Representations of $U_q(\hat{\gl} (n))$ and Quantum Boson-Fermion Correspondence},
Comm. Math. Phys. {\bf 200} (1999), 399--420;
\href{https://arxiv.org/abs/q-alg/9510014}{arXiv:q-alg/9510014}.

\bibitem{DF}
J. Ding, I. B. Frenkel,
{\em Isomorphism of two realizations of quantum affine algebra $U_q(\wht{\gl} (n))$},
Comm. Math. Phys. {\bf 156} (1993), 277--300.

\bibitem{EK3}
P. Etingof, D. Kazhdan,
{\em Quantization of Lie bialgebras, III}, Selecta Math. (N.S.) \textbf{4} (1998), 233--269;
\href{https://arxiv.org/abs/q-alg/9610030}{arXiv:q-alg/9610030}.

\bibitem{EK5}
P. Etingof, D. Kazhdan,
{\em Quantization of Lie bialgebras, V}, Selecta Math. (N.S.) \textbf{6} (2000), 105--130;
\href{http://arxiv.org/abs/math/9808121}{arXiv:math/9808121 [math.QA]}.

\bibitem{FJ}
I. B. Frenkel, N. Jing, 
{\em Vertex representations of quantum affine algebras}, 
Proc. Natl. Acad. Sci. USA, \textbf{85} (1988), 9373--9377.

\bibitem{FLM}
I. Frenkel, J. Lepowsky, A. Meurman,
{\em Vertex operator algebras and the Monster},
Pure and Applied Mathematics, 134. Academic Press, Inc., Boston, MA, 1988.

\bibitem{FR}
 I. B. Frenkel, N. Yu. Reshetikhin, 
{\em Quantum affine algebras and holonomic difference equations}, Comm. Math. Phys. \textbf{146} (1992), 1--60.

\bibitem{G}
M. Gardini,
{\em Quantum vertex algebras}, Ph.D. thesis, Sapienza -- University of Rome, 2018.

\bibitem{I}
K. Iohara,
{\em Bosonic representations of Yangian double $DY_{\hbar}(\mathfrak{g})$ with $\mathfrak{g}=\mathfrak{gl}_N,\mathfrak{sl}_N$},
J. Phys. A \textbf{29} (1996), 4593--4621;
\href{https://arxiv.org/abs/q-alg/9603033}{arXiv:q-alg/9603033}.

\bibitem{J}
M. Jimbo, 
{\em Quantum $R$-matrix for the generalized Toda system}, 
Comm. Math. Phys. \textbf{102} (1986), 537--547.

\bibitem{JKLT}
N. Jing, F. Kong, H.-S. Li, S. Tan,
{\em $(G,\chi_\phi)$-equivariant $\phi$-coordinated quasi modules for nonlocal vertex algebras},
J. Algebra \textbf{570} (2021), 24--74;
\href{https://arxiv.org/abs/2008.05982}{arXiv:2008.05982 [math.QA]}.


\bibitem{JKMY}
N. Jing, S. Ko\v{z}i\'{c}, A. Molev, F. Yang,
{\em Center of the quantum affine vertex algebra in type $A$},
 J. Algebra \textbf{496} (2018), 138--186;
\href{https://arxiv.org/abs/1603.00237}{arXiv:1603.00237 [math.QA]}.

\bibitem{JLM-C}
 N. Jing, M. Liu, A. Molev, 
{\em Isomorphism between the $R$-matrix and Drinfeld presentations of quantum affine algebra: type $C$}, 
J. Math. Phys. \textbf{61} (2020), 031701, 41 pages;
\href{https://arxiv.org/abs/1903.00204}{arXiv:1903.00204 [math.QA]}.

\bibitem{JLM-BD}
 N. Jing, M. Liu, A. Molev, 
{\em Isomorphism between the $R$-matrix and Drinfeld presentations of quantum affine algebra:  types $B$ and $D$}, 
SIGMA Symmetry Integrability Geom. Methods Appl. \textbf{16} (2020), 043, 49 pages;
\href{https://arxiv.org/abs/1911.03496}{arXiv:1911.03496 [math.QA]}.

\bibitem{JLY}
N. Jing, M. Liu, F. Yang,
{\em Double Yangians of classical types and their vertex representations},
J. Math. Phys. \textbf{61} (2020), 051704, (39 pages);
\href{https://arxiv.org/abs/1810.06484}{arXiv:1810.06484 [math.QA]}.

\bibitem{Kac}
V. Kac,
{\em Vertex algebras for beginners}, 
University Lecture Series, 10. American Mathematical Society, Providence, RI, 1997.

\bibitem{Kas}
C. Kassel,
{\em Quantum Groups}, 
Graduate texts in mathematics; vol. \textbf{155}, Springer-Verlag, 1995.

\bibitem{c11}
S. Ko\v{z}i\'{c}, 
{\em Quantum current algebras associated with rational $R$-matrix}, 
 Adv. Math. {\bf 351} (2019), 1072--1104;
\href{https://arxiv.org/abs/1801.03543}{arXiv:1801.03543 [math.QA]}.

\bibitem{Koz1}
S. Ko\v{z}i\'{c}, 
{\em On the quantum affine vertex algebra associated with trigonometric $R$-matrix}, 
Selecta Math. (N.S.) \textbf{27} (2021) 45 (49 pages); 
\href{https://arxiv.org/abs/1908.06517}{arXiv:1908.06517 [math.QA]}.

\bibitem{Li}
H.-S. Li,
{\em $\hbar$-adic quantum vertex algebras and their modules},
Comm. Math. Phys. {\bf 296} (2010), 475--523;
\href{http://arxiv.org/abs/0812.3156}{arXiv:0812.3156 [math.QA]}.

\bibitem{Li1}
H.-S. Li,
{\em $\phi$-Coordinated Quasi-Modules for Quantum Vertex Algebras},
Comm. Math. Phys. {\bf 308} (2011), 703--741;
\href{https://arxiv.org/abs/0906.2710}{arXiv:0906.2710 [math.QA]}.

\bibitem{RS}
N. Yu. Reshetikhin, M. A. Semenov-Tian-Shansky,
{\em Central extensions of quantum current groups}, 
Lett. Math. Phys. {\bf 19} (1990), 133--142.

 
\end{thebibliography}
\end{document}